\magnification=1100
\hfuzz=10 pt
\input psfig.sty

\font\bigbf=cmbx10 at 16pt
\font\medbf=cmbx10 at 13pt
\def\sqr#1#2{\vbox{\hrule height .#2pt
\hbox{\vrule width .#2pt height #1pt \kern #1pt
\vrule width .#2pt}\hrule height .#2pt }}
\def\square{\sqr74}
\def\endproof{\hphantom{MM}\hfill\llap{$\square$}\goodbreak}
\def\i{\item}
\def\c{\centerline}
\def\n{\noindent}

\def\ve{\varepsilon}
\def\n{\noindent}

\def\O{{\cal O}}

\def\C{{\cal C}}
\def\R{I\!\!R}
\def\ov{\overline}

\def\vs{\vskip 2em}
\def\vsk{\vskip 4em}
\def\v{\vskip 1em}
\def\tg{\hbox{tan}}

\null
\c{\bigbf Infinite Horizon Noncooperative Differential Games}
\vs
\c{\it Alberto Bressan$\,^{(*)}$ and  Fabio S. Priuli$\,^{(**)}$}
\v
\c{(*)~Dept.~of Mathematics, Penn State University, University Park
16802, U.S.A.}
\c{bressan@math.psu.edu}
\v
\c{(**)~S.I.S.S.A., Via Beirut 4, Trieste 34014, ITALY.}
\c{priuli@sissa.it}
\vsk
\n{\bf Abstract.} For a non-cooperative differential game, the
value functions of the various players satisfy a
system of Hamilton-Jacobi equations.  In the present paper, we consider
a class of infinite-horizon games with nonlinear costs
exponentially discounted in time.  By the analysis of the value functions,
we establish the existence of Nash equilibrium solutions in feedback form
and provide results and counterexamples on their uniqueness and
stability.
\vsk

\n{\medbf 1 - Introduction}
\v
Problems of optimal control, or zero-sum differential games,
have been the topic of an extensive literature. In both cases, an effective tool
for the analysis of optimal solutions is provided by the
{\it value function}, which  satisfies a
scalar Hamilton-Jacobi equation.  Typically, this first order P.D.E.
is highly non-linear and solutions may not be smooth.
However, thanks to a very effective comparison principle,
the existence and stability of solutions
can be achieved in great generality
by the theory of viscosity solutions, see
[BC] and references therein.

In comparison, much less is known about non-cooperative differential
games.  In a Nash equilibrium solution, the value functions for the
various players now satisfy not a scalar but a system of Hamilton-Jacobi
equations [F].  For this type of nonlinear systems,
no general theorems on the existence or uniqueness of
solutions are yet known.   A major portion of the literature is
concerned with games having linear dynamics and quadratic costs,
see [WSE],[EW],[AFJ] and [PMC].
In this case, solutions are sought among quadratic functions.
This approach effectively reduces the P.D.E.~problem to a finite dimensional
O.D.E..  However, it does not provide insight on the stability
(or instability) of the solutions w.r.t.~small non-linear perturbations.

In [BS1] the first author studied a class of non-cooperative games
with general terminal payoff, in one space dimension.
Relying on recent advances in the
theory of hyperbolic systems of conservation laws, some results on the
existence and stability of Nash equilibrium solutions
could be obtained.  On the other hand, for games in several space
dimensions and also in various one-dimensional cases,
the analysis in [BS2] shows that the corresponding H-J system
is not hyperbolic, hence ill posed.

In the present paper we begin exploring a class of non-cooperative
differential games in infinite time horizon, with exponentially
discounted costs.  In one space dimension,
the corresponding value functions
satisfy a time-independent system of implicit O.D.E's.
Global solutions are sought within a class of absolutely
continuous functions,
imposing certain growth conditions as $|x|\to\infty$, and suitable
admissibility conditions at points where the gradient $u_x$
has a jump.

The dynamics of our system is very elementary, and the
cost functions that we consider are small perturbations of linear ones.
However, already in this simple setting we find cases where the problem has
unique solution, and cases where infinitely many solutions exist.
This provides a glimpse of the extreme complexity of the problem,
for general non-cooperative $N$-player games with non-linear cost functions.

The plan of the paper is as follows. In Section 2 we describe the
differential game, introducing
the basic notations and definitions.  In Section 3 we prove that, from
an admissible solution to the O.D.E.~for the value function,
one can always recover a Nash equilibrium solution to the differential game.
The relevance of our admissibility conditions is then highlighted
by two examples.

The existence and uniqueness of global admissible solutions to the
H-J system for the value functions is then studied in Sections 3 and 4.
We first consider the cooperative case, where both players wish to
move the state of the system in the same direction. In the case with
terminal payoff, this situation was leading to a
well-posed hyperbolic Cauchy problem [BS1].  As expected, in
the infinite-horizon case
we still obtain an existence and uniqueness result.
Subsequently, we consider the case of conflicting interests, where
the players wish to steer the system in opposite directions.
In the case with
terminal payoff, this situation leads to an
ill-posed Cauchy problem, as shown in [BS2].
Somewhat surprisingly, we find that
the corresponding infinite-horizon case can have unique or multiple solutions,
depending on the values of certain parameters.
\vsk
\n{\medbf 2 - Basic definitions}
\v
Consider an
$m$-persons non-cooperative differential game, with dynamics
$$\dot x=\sum_{i=1}^m f_i(x,\alpha_i),\qquad\qquad
\alpha_i(t)\in A_i\,,\qquad\qquad x\in\R^n.\eqno(2.1)$$
Here $t\mapsto\alpha_i(t)$
is the control chosen by the $i$-th player, within a set
of admissible control values $A_i\subseteq \R^k$.
We will study the {\bf discounted, infinite horizon problem}, where
the game takes place on an infinite interval of time $[0,\,\infty[\,$,
and each player has only a running cost, discounted exponentially
in time.
More precisely,
for a given initial data
$$x(0)=y\in\R^n\,,\eqno(2.2)$$
the goal of the $i$-th player is to minimize the functional
$$J_i(\alpha)\doteq \int_0^\infty
e^{-t}\,\psi_i\big(x(t),\,\alpha_i(t)\big)\,dt\,,\eqno(2.3)$$
where $t\mapsto x(t)$ is the trajectory of (2.1).
By definition, an $m$-tuple of feedback strategies
$\alpha_i=\alpha_i^*(x)$, $i=1,\ldots,m$,
represents a {\it Nash
non-cooperative equilibrium solution}
for the differential game (2.1)-(2.2) if the following holds.
For every $i\in \{1,\ldots,m\}$,
the feedback control $\alpha_i=\alpha^*(x)$
provides a solution to the
the optimal control problem for the
$i$-th player,
$$\min_{\alpha(\cdot)} ~J_i(\alpha)\,,\eqno(2.4)$$
where the dynamics of the system is
$$\dot x=f_i(x,\alpha_i)+\sum_{j\not= i} f_j(x,\alpha^*_j(x)),
\qquad\qquad \alpha_i(t)\in A_i\,.\eqno(2.5)$$
More precisely, we require that, for every initial data
$y\in\R$, the Cauchy problem
$$\dot x=\sum_{j=1}^m f_j\big(x,\alpha_j^*(x)\big)\,,\qquad
\qquad x(0)=y\,,
\eqno(2.6)$$
should have at least one Caratheodory solution
$t\mapsto x(t)$, defined for all $t\in [0,\infty[\,$.
Moreover, for every such solution and each $i=1,\ldots,m$,
the cost to the $i$-th player
should provide the minimum for the optimal control problem
(2.4)-(2.5).
We recall that a Caratheodory solution is an absolutely
continuous function $t\mapsto x(t)$
which satisfies the differential equation
in (2.6) at almost every $t> 0$.

Nash equilibrium solutions in feedback form can be
obtained by studying a related system of P.D.E's.
Assume that a value function $u(y)=(u_1,\ldots,u_n)(y)$ exists,
so that $u_i(y)$ represents the cost for the $i$-th player
when the initial state of the system is $x(0)=y$ and
the strategies $\alpha_1^*,\ldots,\alpha_m^*$ are implemented.
By the theory of optimal control, see for example [BC],
on regions where $u$ is smooth,
each component $u_i$ should provide a solution
to the corresponding scalar Hamilton-Jacobi-Bellman equation.
The vector function $u$ thus satisfies the
stationary system of equations
$$u_i(x)=H_i(x,\,\nabla u_1,\ldots,\nabla u_m),\eqno(2.7)$$
where the Hamiltonian functions $H_i$ are defined as follows.
For each $p_j\in\R^n$, assume that there exists an optimal control
value $\alpha_j^*(x,p_j)$ such that
$$p_j\cdot f_j\big(x,\,\alpha_j^*(x,p_j)\big)+
\psi_j\big(x,\,\alpha_j^*(x,p_j)\big)=\min_{a\in A_j}
\,\big\{p_j\cdot f_j(x,a)+
\psi_j(x,a)\big\}\,.\eqno(2.8)$$
Then
$$H_i(x,\,p_1,\ldots,p_m)\doteq p_i\cdot\sum_{j=1}^m f_j
\big(x,\,\alpha_j^*(x,p_j)\big)+\psi_i\big(x,\,\alpha_i^*(x,p_i)\big)\,.
\eqno(2.9)$$

A rich literature is currently available on optimal control
problems and on viscosity solutions
to the corresponding scalar H-J equations.
However, little is yet known about non-cooperative differential games,
apart from the linear-quadratic case.
In this paper we begin a study of this class of differential games,
with two players in one space dimension.  Our main interest is in the
existence, uniqueness and stability of Nash equilibrium solutions
in feedback form.

When $x$ is a scalar variable, (2.7) reduces to a system
of implicit O.D.E's:
$$u_i=H_i(x, u_1',\ldots, u_m')\,.\eqno(2.10)$$
In general, this system will have infinitely many solutions.
To single out a (hopefully unique) admissible solution,
corresponding to a Nash equilibrium for the differential game,
additional requirements must be imposed. These are of two types:
\v
\i{(i)} Asymptotic growth conditions as $|x|\to\infty$.
\v
\i{(ii)}
Jump conditions, at points where the derivative $u'$ is discontinuous.
\v
\n
To fix the ideas, consider a game with the simple dynamics
$$\dot x(t)= \alpha_1(t)+\cdots +\alpha_m(t)\,,\eqno(2.11)$$
and with cost functionals of the form
$$J_i(\alpha)\doteq\int_0^\infty e^{-t} \bigg[ h_i\big(x(t))
+k_i\big(x(t)\big)\,
{\alpha_i^2(t)\over 2}
\bigg]\,dt\,.\eqno(2.12)$$
We shall assume that the functions $h_i, k_i$ are smooth and satisfy
$$\big|h_i'(x)\big|\leq C\,,\qquad\qquad
{1\over C}\leq k_i(x)\leq C\,,\eqno(2.13)$$
for some constant $C>0$.
Notice that in this case (2.8)
yields $\alpha_i^*=-p_i/k_i$, hence (2.10) becomes
$$u_i=\left({u_i'\over 2k_i(x)}-\sum_{j=1}^m {u_j'\over k_j(x)}
\right)u_i'+h_i(x)\,.\eqno(2.14)$$
For a solution to the system of H-J equations (2.14),
a natural set of admissibility conditions is formulated below.
\v
\n{\bf Definition 1.}
A function
$u:\R\mapsto\R^m$ is called an {\it admissible solution} to the
implicit system of O.D.E's (2.14) if the following holds.
\v
\n(A1) ~$u$ is absolutely continuous.
Its derivative $u'$ satisfies the
equations (2.14) at a.e.~point $x\in\R$.
\v
\n(A2) ~$u$ has sublinear growth at infinity. Namely,
there exists a constant $C$ such that, for all $x\in\R$,
$$\big|u(x)\big|\leq C\,\big(1+|x|\big)\,.\eqno(2.15)$$
\v
\n(A3) At every point $y\in \R$, the derivative $u'$ admits right and
left limits $u'(y+)$, $u'(y-)$.  At points where $u'$ is discontinuous,
these limits satisfy the
admissibility conditions
$$\sum_{i=1}^m {u_i'(y+)\over k_i(y)}~\leq~ 0~\leq~
\sum_{i=1}^m {u_i'(y-)\over k_i(y)}\,.\eqno(2.16)$$
\v
Because of the assumption (2.13), the cost functions
$h_i$ are globally Lipschitz continuous. It is thus natural to
require that the value functions $u_i$ be absolutely
continuous, with sub-linear growth
as $x\to\pm\infty$.   Call $p_i^\pm\doteq u_i'(y\pm)$.
By the equations (2.14) and the continuity of the functions $u_i, h_i, k_i$,
one obtains the identities
$${(p_i^+)^2\over 2k_i(y)}+\sum_{j\not= i}{p_i^+p_j^+\over
k_j(y)}~=~
{(p_i^-)^2\over 2k_i(y)}+\sum_{j\not= i}{p_i^-p_j^-\over
k_j(y)}\qquad\qquad i=1,\ldots,m\,.\eqno(2.17)$$
Recalling that the feedback controls are $\alpha_i^*=-u_i'/k_i$,
the condition (2.16) now becomes clear: it states that
$\dot x(y-)\leq 0\leq \dot x(y+)$, i.e.,~trajectories
should move away from a point of discontinuity.

Notice that all of the above conditions
are satisfied at a point $y$ such that
$$ \sum_j {u_j'(y+)\over k_j(y)}\leq 0\,,
\qquad\qquad u'_i(y+)+u'_i(y-)=0\qquad i=1,\ldots,m\, .\eqno(2.18)$$

By (A1), the derivatives $p_i=u_i'$ are defined at
a.~e.~point $x\in\R$. The optimal feedback controls
$\alpha_i^*=-p_i/k_i$ are thus defined almost everywhere.
We can use the further assumption (A3) and extend these functions
to the whole real line by taking  limits from the right:
$$\alpha^*(x)\doteq  -{u_i'(x+)\over k_i(x)}\,.\eqno(2.19)$$
In this way, all feedback control functions will be right-continuous.

The system of implicit differential equations (2.14) is highly nonlinear
and difficult to study in full generality.
In this paper we initiate the analysis by looking at some
significant cases. Our main results can be roughly
summarized as follows:
\v
\i{(i)} If $u=(u_1,\ldots u_m)$ provides an admissible solution
to the system of Hamilton-Jacobi equations (2.14),
then the feedback strategies (2.19)
provide a Nash equilibrium solution to the differential game
(2.11)-(2.12).
\v
\i{(ii)} For games with two players, if the cost functions
$h_1, h_2$ are both monotone increasing
(or both monotone decreasing), then (2.14) has
a unique admissible solution.
\v
\i{(iii)} Still in the case of two players, one can give examples where
the derivatives of the cost functions
satisfy $h_i'+h_2'=0$ and infinitely many admissible solutions
of (2.14) are found. On the other hand, if the sum $h_1'+h_2'$ remains
bounded away from zero, then under suitable assumptions the system
(2.14) has a unique admissible solution.
\vsk

\n{\medbf 3 - Solutions of the differential game}
\v
In this section we prove that admissible solutions to the H-J equations
yield a solution to the differential game. Moreover, we give a couple
of examples showing the relevance of the assumptions (A2) and (A3).

\v
\n{\bf Theorem 1.} {\it Consider the differential game (2.11)-(2.12),
with the assumptions (2.13).
Let $u:\R\mapsto \R^m$ be an admissible solution to the
systems of H-J equations (2.14), so that the conditions (A1)--(A3)
hold.  Then the controls (2.19)
provide a Nash equilibrium solution in feedback form.}
\v
\n{\bf Proof.}
The theorem will be proved in several steps.
\v
\n{\bf 1.} First of all, setting
$$g(x)~\doteq~ \sum_i \alpha_i^*(x)~=~- \sum_i {u_i'(x)\over k_i(x)}\,,
\eqno(3.1)$$
we need to prove that the Cauchy problem
$$\dot x(t)=g\big(x(t)\big)\,,
\qquad\qquad x(0)=y\,,\eqno(3.2)$$
has a globally defined solution, for every initial data $y\in \R$.
This is not entirely obvious, because the function $g$ may be discontinuous.
We start by proving the local existence of solutions.
\v
\n CASE 1: $g(y)=0$. In this trivial case $x(t)\equiv y$ is the required solution.
\v
\n CASE 2: $g(y)>0$. By right continuity, we then have $g(x)>0$ for
$x\in [y,~y+\delta]$, for some $\delta>0$.  This implies the existence
of a (unique) strictly increasing solution $x:[0,\ve]\mapsto
\R$, for some $\ve>0$.
\v
\n CASE 3: $g(y)<0$.  By the admissibility conditions (2.16),
this implies that $g$ is continuous and negative in a neighborhood
of $y$.  Therefore the Cauchy problem (3.2) admits
a (unique) strictly decreasing solution $x:[0,\ve]\mapsto
\R$, for some $\ve>0$.
\v
\n{\bf 2.} Next, we prove that the local solution can be extended
to all positive times.   For this purpose, we need to rule out the
possibility that
$\big|x(t)\big|\to\infty$ in finite time.
We first observe that each trajectory is monotone, i.e., either non-increasing,
or non-decreasing, for $t\in [0,\infty[\,$.
To fix the ideas,  let $t\mapsto x(t)$ be strictly increasing,
with $x(t)\to\infty$ as $t\to T-\,$.   A contradiction is now
obtained as follows.  For each $\tau>0$, using (2.14) we compute
$$\eqalign{\sum_i & u_i\big(x(\tau)\big)-\sum_i u_i\big(x(0)\big)
=\int_0^\tau \left\{ {d\over dt} \sum_i u_i\big(x(t)\big)\right\}\,dt\cr
&=\int_0^\tau -\left\{\sum_i u_i'\big(x(t)\big)\cdot
\sum_j {u_j'\big(x(t)\big)\over k_j\big(x(t)\big)}\right\}\,dt\cr
&=\int_0^\tau \sum_i \left\{ u_i\big(x(t)\big)
-{\big|u_i'(x(t))\big|^2\over 2 k_i(x(t))}- h_i\big(x(t)\big)\right\}
\,dt\cr}\eqno(3.3)$$
By assumptions, the functions $u_i$ and $h_i$ have sub-linear growth.
Moreover, each $k_i$ is uniformly positive and bounded above.
Using the elementary inequality
$$\big|x(\tau)-x(0)\big|\leq \int_0^\tau 1\cdot \big|\dot x(t)\big|\,dt
\leq \left(\int_0^\tau 1\,dt\right)^{1/2}\cdot\left(\int_0^\tau
\big|\dot x(t)\big|^2\,dt\right)^{1/2},$$
from
(3.3) we thus obtain
$$\eqalign{ {\big|x(\tau)-x(0)\big|^2\over\tau}&
\leq \int_0^\tau \big|\dot x(t)\big|^2\,dt
~\leq 4C \int_0^\tau \sum_i{\big|u_i'(x(t))\big|^2\over 2 k_i(x(t))}\,dt\cr
&\leq 4C \left\{~\left|\sum_i  u_i\big(x(\tau)\big)-\sum_i u_i\big(x(0)\big)\right|
+ \int_0^\tau \sum_i \big| u_i(x(t))\big|\,dt
+ \int_0^\tau\sum_i \big| h_i(x(t))\big|
\,dt\right\}\cr
&\leq C_0\, (1+\tau)\Big\{ 2+\big|x(\tau)\big|+\big|x(0)\big|\Big\}\,,}
$$
for some constant $C_0$.  Therefore, either $\big| x(\tau)\big|\leq 2+3
\big|x(0)\big|$,
or else
$$\big|x(\tau)\big|\leq \big|x(0)\big|+2\tau\cdot C_0\,(1+\tau)\,.\eqno(3.4)$$
In any case, blow-up cannot occur at any finite time $T$.
\v
\n{\bf 3.} To complete the proof,
for each fixed $i\in\{1,\ldots,m\}$,
we have to show that the feedback $\alpha_i^*$ in (2.19)
provides solution to the  optimal control problem for the $i$-th player:
$$\min_{\alpha_i(\cdot)} \int_0^\infty e^{-t} \bigg[ h_i\big(x(t))
+k_i\big(x(t)\big)\,
{\alpha_i^2(t)\over 2}
\bigg]\,dt\,,\eqno(3.5)$$
where the system has dynamics
$$\dot x=\alpha_i+\sum_{j\not= i} \alpha_j^*(x)\,.\eqno(3.6)$$
Given an initial state $x(0)=y$,
by the assumptions on $u$ it follows that the feedback strategy
$\alpha_i=\alpha_i^*(x)$ achieves a total cost given by
$u_i(y)$.  Now consider any absolutely continuous trajectory
$t\mapsto x(t)$, with $x(0)=y$. Of course, this corresponds to
the control
$$\alpha_i(t)\doteq \dot x(t)-\sum_{j\not= i} \alpha_j^*(x)\eqno(3.7)$$
implemented by the $i$-th player.
We claim that the corresponding cost satisfies
$$
\int_0^\infty e^{-t} \bigg[ h_i\big(x(t)\big)
+{k_i\over 2}\Big(\dot x(t)-\sum_{j\not= i}
\alpha_j^*\big(x(t)\big)\Big)^2
\bigg]\,dt~\geq ~u_i(y)\,.\eqno(3.8)$$
To prove (3.8), we first observe that (3.4) implies
$$\lim_{t\to\infty} e^{-t} u_i\big(x(t)\big)=0
\qquad\qquad i=1,\ldots,n\,.$$
Hence
$$u_i(y)=u_i\big(x(0)\big)=-\int_0^\infty {d\over dt}
\bigg[e^{-t} u_i\big(x(t)\big)\bigg]\,dt\,.$$
The inequality (3.8) can now be established by checking that
$$e^{-t}\bigg[ h_i\big(x(t)\big)
+{k_i\over 2}\Big(\dot x(t)-\sum_{j\not= i}
\alpha_j^*\big(x(t)\big)\Big)^2
\bigg]\geq e^{-t}u_i\big(x(t)\big)-e^{-t}
u_i'\big(x(t)\big)\cdot \dot x(t)\,.\eqno(3.9)$$
Equivalently, letting $\alpha_i$ be as in (3.7),
$$u_i\leq\left( \alpha_i-\sum_{j\not= i} {u_j'\over k_j}\right)
u_i'+{k_i\over 2}\alpha_i^2+h_i\,.$$
This is clearly true because, by (2.8),
$$u_i(x)=\min_{a} \left\{{k_i\over 2} a^2+a u_i'
-\sum_{j\not= i}{u_j'u_i'\over k_j} +h_i(x)\right\}\,.$$
\endproof

\vsk
We now give two examples showing that,
if the growth assumptions (2.15) or if the jump conditions
(2.16) are not satisfied, then the feedbacks (2.19) may not
provide a Nash equilibrium solution. This situation is well known already
in the context of control problem.
\v
\n{\bf Examples.} Consider the game for two players, with
dynamics
$$\dot x=\alpha_1+\alpha_2\,,\eqno(3.10)$$
and cost functionals
$$J_i=\int_0^\infty e^{-t}\cdot {\alpha_i^2(t)\over 2}
\,dt\,.$$
In this case, if $u_i'=p_i$, the optimal control for the $i$-th player is
$$\alpha_i^*(p_i)~=~\hbox{arg}~\min_\omega\Big\{
p_i\,\omega +{\omega^2\over 2}\Big\}~=~- p_i\,.$$
The system of H-J takes the simple form
$$\left\{ \eqalign{u_1&= -\left( {u_1'\over 2}+u_2'\right)u_1'\,,\cr
u_2&= -\left( u_1'+{u_2'\over 2}\right)u_2'\,.\cr}\right.\eqno(3.11)$$
The obvious admissible solution is $u_1\equiv u_2\equiv 0$,
corresponding to identically zero controls, and zero cost.
We now observe that
the functions
$$u_1(x)=\cases{ 0\quad &if\quad $|x|\geq 1\,$,\cr
-{1 \over 2}\big(1-|x|\big)^2\quad &if\quad $|x|< 1\,$,\cr
}\qquad\qquad u_2(x)=0\,,$$
provide a solution to (3.11), which is not admissible
because the conditions (2.16) fail at $x=0$.

Next, the functions
$$u_1(x)=-{1 \over 2}x^2,
\qquad\qquad u_2(x)=0\,,$$
provide yet another another solution,
which does not satisfy the growth conditions
(2.15).

In the above two cases, the corresponding feedbacks $\alpha_i^*(x)=-u_i'(x)$
do not yield
a solution to the differential game.
\vsk
\n{\medbf 4 - Cooperative situations}
\v
We consider here a game for two players, with dynamics
$$\dot x=\alpha_1+\alpha_2\,,\eqno(4.1)$$
and cost functionals of the form
$$J_i(\alpha)\doteq\int_0^\infty e^{-t}
\bigg[ h_i\big(x(t)) +
{\alpha_i^2(t)\over 2}
\bigg]\,dt\,.\eqno(4.2)$$
Notice that, for any positive constants $k_1,k_2,\lambda$, the more general
case
$$J_i(\alpha)\doteq\int_0^\infty e^{-\lambda t}
\bigg[ \tilde h_i\big(x(t)) +
{\alpha_i^2(t)\over 2 k_i}
\bigg]\,dt$$
can be reduced to (4.2) by a linear change of variables.

The
system of H-J equations
for the value functions now takes the form
$$\left\{\eqalign{
u_1(x)&=h_1(x)-u_1'u_2'-(u_1')^2/2\,,\cr
u_2(x)&=h_2(x)-u_1'u_2'-(u_2')^2/2\,,\cr
}\right.\eqno(4.3)$$
and the optimal feedback controls are given by
$$\alpha_i^*(x)=-u_i'(x)\,.\eqno(4.4)$$
Differentiating (4.3) and setting $p_i=u_i'$
one obtains
the system
$$\left\{ \eqalign{ h_1'-p_1&= (p_1+p_2)p_1'+p_1p_2'\,,\cr
h_2'-p_2&= p_2p_1'+(p_1+p_2)p_2'\,.
\cr}\right.
\eqno(4.5)$$
Set
$$\Lambda(p)\doteq\pmatrix{ p_1+p_2 & p_1\cr
p_2 & p_1+p_2\cr}\,, \qquad\qquad\Delta(p)\doteq\det\,\Lambda(p)\,.$$
{}From (4.5) we deduce
$$\left\{
\eqalign{p'_1
&=\Delta(p)^{-1}\big[- p_1^2 +(h'_1-h'_2) p_1+h'_1p_2\big]\,,\cr
p'_2
&=\Delta(p)^{-1}\big[-p_2^2 +(h'_2-h'_1) p_2+h'_2p_1
\big]\,.\cr}
\right.\eqno(4.6)$$
Notice that
$$ {1\over 2}\,(p_1^2+p_2^2)\leq \Delta(p)\leq 2(p_1^2+p_2^2)\,.\eqno(4.7)$$

\n In particular, $\Delta(p)>0$ for all $p=(p_1,p_2)\not= (0,0)$. Up
to a rescaling of the independent variable, we can thus study the
equivalent system
$$\left\{
\eqalign{p_1'
&=(h'_1-h'_2) p_1+h'_1p_2-p_1^2\,,\cr
p_2'
&=(h'_2-h'_1) p_2+h'_2p_1-p_2^2
\,.\cr}
\right.\eqno(4.8)$$
For piecewise smooth solutions, jumps are only allowed from
any point $(p_1^-, p_2^-)$ with
$$p_1^-+p_2^-\geq 0\eqno(4.9)$$
to the symmetric point
$$(p_1^+,p_2^+)= (-p_1^-,~- p_2^-)\,.\eqno(4.10)$$

\v
\n{\bf Theorem 2.} {\it Let the cost functions $h_1, h_2$
be smooth, and assume that their derivatives satisfy
$${1\over C}\leq h_i'(x)\leq C\eqno(4.11)$$
for some constant $C>1$ and all $x\in\R$.
Then the system (4.3)
has an admissible solution. The corresponding
functions
$\alpha_i^*$ in (4.4) provide a Nash equilibrium solution
to the non-cooperative game.}
\v
\n{\bf Proof.}
Write the O.D.E.~(4.6) in the more compact form
$${dp\over dx} = f(p)\,.\eqno(4.12)$$
To show the existence of at least one admissible solution of (4.3),
for
every $\nu\geq 1$ let $p^{(\nu)}:\,[-\nu,\infty[\,\to\R^2$
be the solution of the
Cauchy problem
$${dp^{(\nu)}\over dx} = f\big(p^{(\nu)}\big)\,,\qquad\qquad
p^{(\nu)}(-\nu)=(1,1).\eqno(4.13)$$
It is easy to check that the polygon
$$\Gamma\doteq \Big\{ (p_1,p_2)\,; ~~p_1,p_2\in [0,2C]\,,~~p_1+p_2\geq 1/2C
\Big\}$$
is positively invariant for the flow of (4.6).
Hence $p^{(\nu)}(x)\in\Gamma$ for all $\nu\geq 1$ and $x\geq -\nu$.
We can extend each function $p^{(\nu)}$
to the whole real line by setting
$$p^{(\nu)}(x)=(1,1)\qquad\qquad\hbox{for}~~x<-
\nu\,.$$
By uniform boundedness and equicontinuity, the sequence $p^{(\nu)}$
admits a subsequence converging to a uniformly continuous function
$p:\R\mapsto\Gamma$.
Clearly this limit function provides a continuous,
globally bounded solution of (4.6).
We then define the controls $\alpha_i^*(x)\doteq -p_i(x)$
and the
cost functions
$$u_i(y)\doteq
\int_0^\infty e^{-t}
\bigg[ h_i\big(x(t,y)) +
{1\over 2} \big(\alpha_i^*\big(x(t,y)\big)\Big)^2
\bigg]\,dt\,,\eqno(4.14)$$
where $t\mapsto x(t,y)$ denotes the solution to the
Cauchy problem
$$\dot x=\alpha_1^*(x)+\alpha_2^*(x)\,,\qquad\qquad x(0)=y\,.
\eqno(4.15)$$
This function provides a globally Lipschitz, smooth solution of
the system (4.3).
\endproof
\v
In the case where the oscillation of the derivatives $h_i'$
is sufficiently small, we can also prove the uniqueness of the
Nash feedback solution.

\v
\n{\bf Theorem 3.} {\it Let the cost functions be smooth, with derivatives
satisfying (4.11), for some constant $C$.
Assume that the oscillation of their derivatives satisfies
$$\sup_{x,y\in\R} \big|h_i'(x)-h_i'(y)\big|\leq\delta \qquad i=1,2
\eqno(4.16)$$
for some $\delta>0$ sufficiently small (depending only on $C$).
Then the admissible solution of
the system (4.3) is unique.
}
\v
Before giving details of the proof, we sketch the main ideas.
In the case of linear cost functions, where $h_i(x)=\kappa_i \,x$,
$h_i'\equiv \kappa_i$,
the phase portrait of the planar O.D.E. (4.8) is depicted in
Figure 1. We observe that

- Unbounded trajectories of (4.8), with $\big|p(s)\big|\to \infty$ as
$s\to \bar s$, correspond to solutions $p=p(x)$ of (4.6) with
$\big|p(x)\big|\to \infty$, $\big|p'(x)\big|\to \infty$
as $x\to\pm\infty $.   Indeed, because of the rescaling
(4.7), as the parameter $s$ approaches s finite limit $\bar s$,
we have  $|x|\to\infty$.  This yields a solution $u(x)=\int_*^x p(x)\,dx$
which does not satisfy the growth restrictions (2.15).

- The heteroclinic orbit, joining the origin
with the point $(\kappa_1,\kappa_2)$, corresponds to a trajectory of
(4.16) defined on a half line, say $[\bar x, \infty[\,.$
To prolong this solution for $x<\bar x$ one needs a trajectory
of (4.16) which approaches the origin as $s\to\infty$. But the two available
solutions are both unbounded, hence not acceptable.

- Finally, one must examine
solutions whose gradient has one or more jumps, from
a point $P=(p_1,p_2)$ with $p_1+p_2\geq 0$ to its symmetric
point $-P=(-p_1,\,-p_2)$. However, a direct inspection shows that,
even allowing these jumps, one still cannot construct
any new globally bounded trajectory.

In the end, in linear case,
one finds that the only admissible solution is $(p_1,p_2)\equiv
(\kappa_1,\kappa_2)$.  A perturbative argument shows that this
conclusion remains valid if a small  $\C^1$ perturbation is added to the
cost functions.
\vs
\n{\bf Proof of Theorem 3.} {\it First Step.} We begin with the case
$h'_i(x)\equiv\kappa_i$ and assume, without any loss of generality,
that $\kappa_1\leq\kappa_2$.


Let $\tilde p$ be a smooth solution
of (4.8), as shown in Figure 1.

\midinsert \vskip 10pt
\centerline{\hbox{\psfig{figure=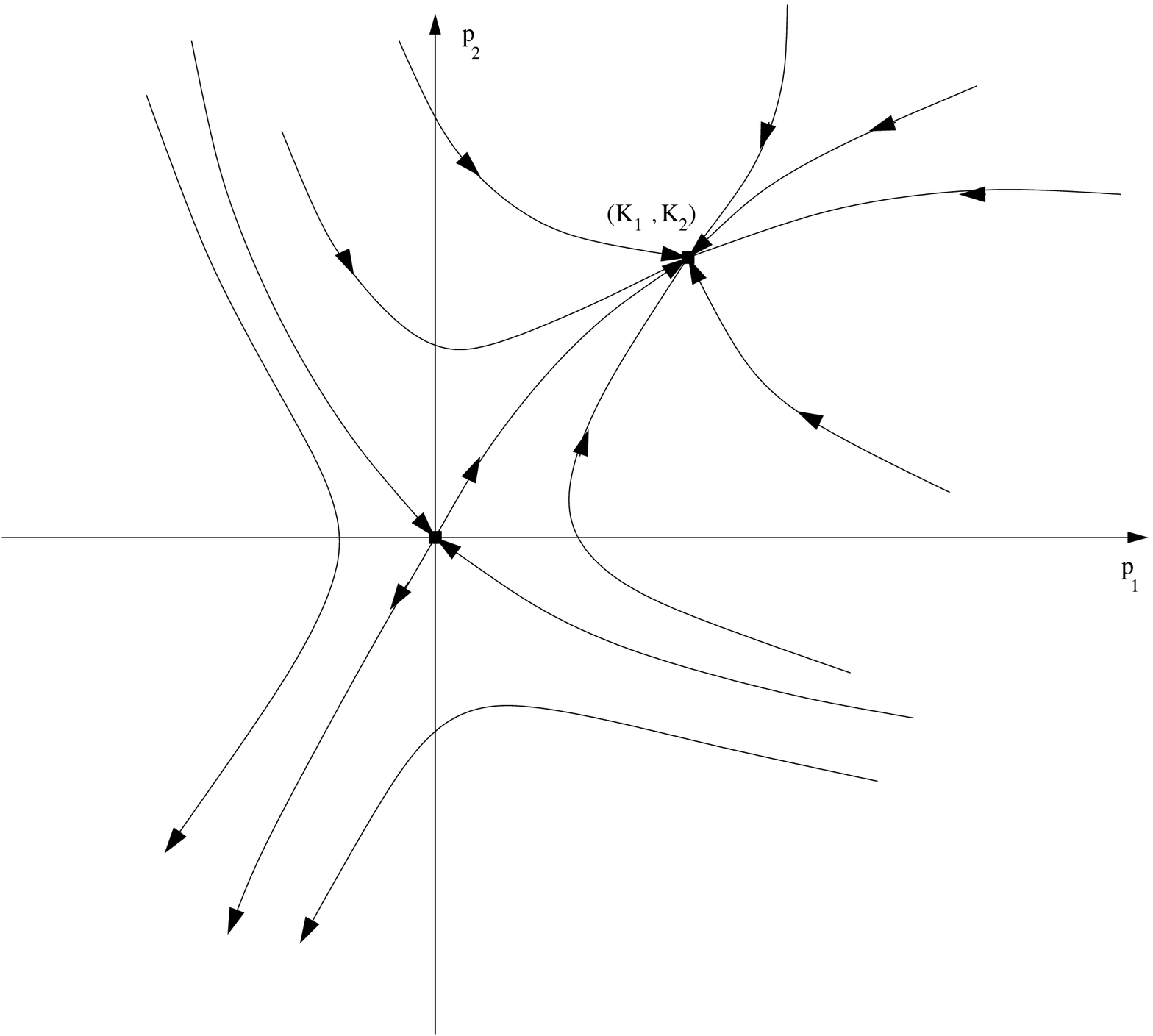,width=12cm}}}
\centerline{\hbox{Figure 1}} \vskip 10pt
\endinsert

\n We observe that
the following facts hold (see Figure~2):

\v \n {\bf 1.}   Both sets
$A=\{ (p_1,p_2)\neq (0,0)\,:~p_1\leq 0\,,~p_2\leq 0\}$ and
$\{ (p_1,p_2)\neq (0,0)\,: ~p_1\geq 0\,,~p_2\geq 0\}$ are
positively invariant for the flow of (4.8) and both
$B=\{ (p_1,p_2)\,:~p_1>0\,,~p_2< 0\}$ and
$C=\{ (p_1,p_2)\,:~p_1<0\,,~p_2> 0\}$ are negatively invariant.

\v \n {\bf 2.}   If $\tilde p (s_o)\in A=\{ (p_1,p_2)\neq
(0,0)\,: ~p_1\leq 0\,,~p_2\leq 0\}$ for some $s_o$, then
$|\tilde p|\to +\infty$ as $s\to +\infty$. Indeed, since
$$
{d\over ds} (\tilde p_1+\tilde p_2)=
-\tilde p_1^2-\tilde p_2^2+\kappa_1 \tilde p_1+\kappa_2
\tilde p_2 \leq -{1 \over 2} (\tilde p_1+\tilde p_2)^2 <0\,,
$$
we can assume there exist $\bar s\geq s_o$ and $\ve>0$ such
that $\tilde p_1(\bar s)+\tilde p_2(\bar s)<-\ve$.
Moreover, the following holds for any $\sigma>\bar s$:
$$
{d\over ds} (\tilde p_1+\tilde p_2) (\sigma)\leq -{1 \over 2}
(\tilde p_1 (\sigma)+\tilde p_2 (\sigma))^2 < -{1 \over 2}(\tilde p_1
(\bar s)+\tilde p_2 (\bar s))(\tilde p_1 (\sigma)+\tilde p_2 (\sigma))<
{\ve \over 2}(\tilde p_1 (\sigma)+\tilde p_2 (\sigma))\,.
$$
Hence, an integration yields $(\tilde p_1+\tilde p_2)(s)\leq -\eta
e^{{\ve \over 2} s}$ for $s>\bar s$ (and $\eta>0$) and
$(\tilde p_1+\tilde p_2)\to -\infty$ as $s\to +\infty$.

\v \n {\bf 3.}   If $\tilde p(s_o)\in B=\{ (p_1,p_2)\,:
~p_1>0\,,~p_2< 0\}$ for some $s_o$, then $|\tilde p|\to
+\infty$ as $s\to -\infty$. Indeed, let $\ve>0$ such that
$\tilde p_1(s_o)>\ve$. Since
$$
{d\over ds} \tilde p_1 =-\tilde p_1^2+(\kappa_1-\kappa_2)\tilde p_1+\kappa_1
\tilde p_2 \leq - (\tilde p_1+\kappa_2-\kappa_1)
\tilde p_1<0\,,
$$
it is sufficient to observe that, for $\sigma< s_o$,
$$
{d\over ds} \tilde p_1 (\sigma) <
-(\tilde p_1(s_o)+\kappa_2-\kappa_1) \tilde p_1(\sigma) \leq -
(\ve+\kappa_2-\kappa_1) \tilde p_1(\sigma)\,.
$$
Hence, an integration yields $\tilde p_1(s)\geq \eta
e^{- (\ve+\kappa_2-\kappa_1) s}$ for $s<s_o$ (and $\eta>0$) and
$\tilde p_1\to +\infty$ as $s\to -\infty$.

\midinsert \vskip 10pt
\centerline{\hbox{\psfig{figure=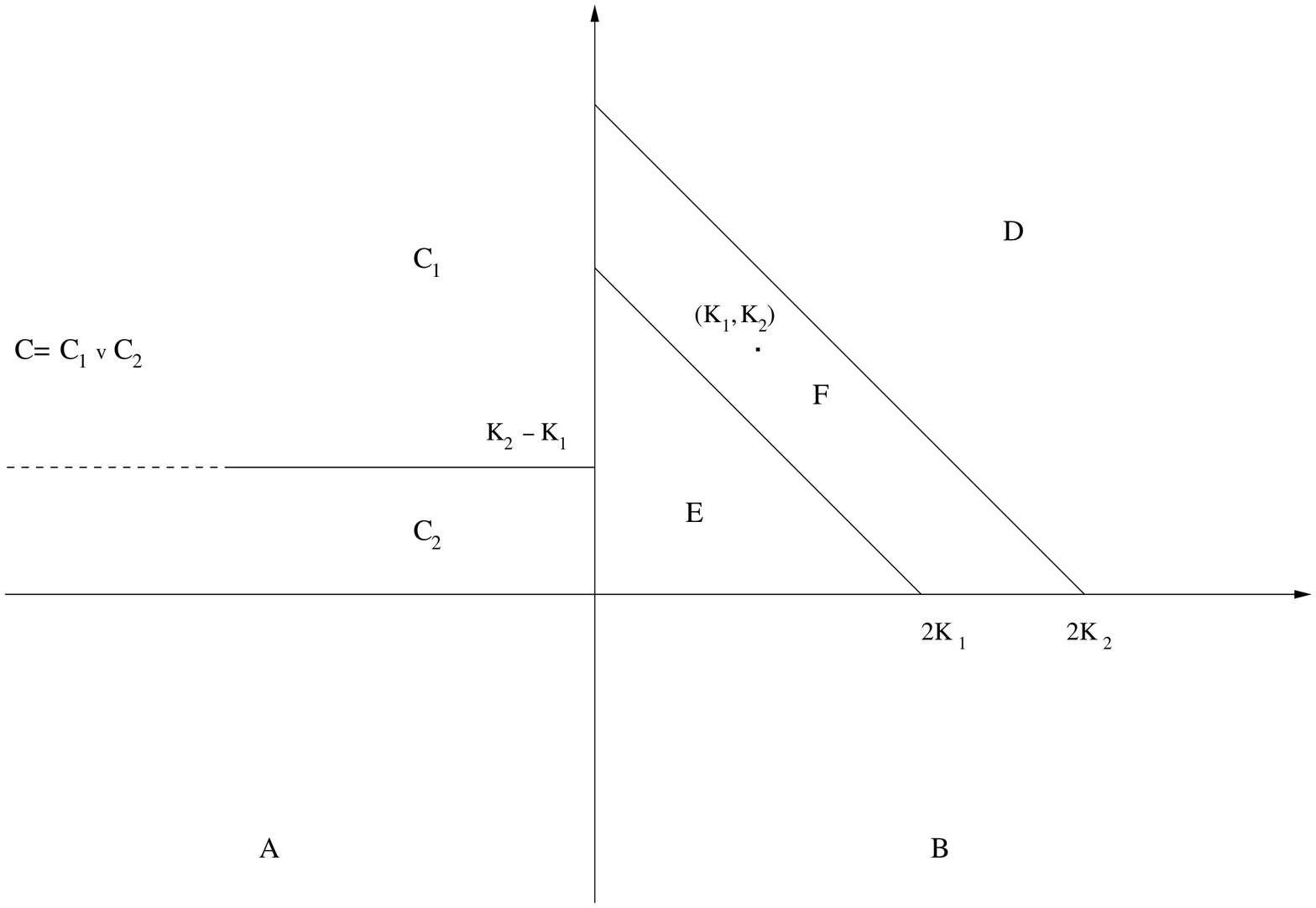,width=10cm}}}
\centerline{\hbox{Figure 2}} \vskip 10pt
\endinsert

\v \n {\bf 4.}   If $\tilde p(s_o)\in C_1=\{ (p_1,p_2)\,:
~p_1< 0\,,~p_2>\kappa_2-\kappa_1\}$ for some $s_o$, then
$|\tilde p|\to +\infty$ as $s\to -\infty$. Here the argument is
exactly the same as in the previous case with $\tilde p_2$ in place of
$\tilde p_1$.

\v \n {\bf 5.}   If $\tilde p(s_o)\in C_2=\{ (p_1,p_2)\,:
~p_1< 0\,,~0<p_2\leq\kappa_2-\kappa_1\}$ then there exists $\bar s< s_o$ such
that  $\tilde p(\bar s)$ is in $C_1$ as in case {\bf 4.} above. Indeed
there could be only two situations.

\n If $-\tilde p_1(s_o)^2+(\kappa_1-\kappa_2)\tilde p_1(s_o)+\kappa_1
\tilde p_2(s_o)\geq 0$, then, by negative invariance, $\tilde p$ could
only have reached this region from $C_1$, hence there
exists $\bar s\leq s_o$ such that $\tilde p(\bar s)$ is as in case
{\bf 4} above. Otherwise, using again negative invariance and the fact that
there are no equilibria in $C_2$, either there exists
$\bar s\leq s_o$ such that $\tilde p(\bar s)$ is in case {\bf 4} above,
or there exists $s_1< s_o$ such that $-\tilde p_1(s_1)^2+
(\kappa_1-\kappa_2)\tilde p_1(s_1)+\kappa_1\tilde p_2(s_1)\geq 0$ and then,
by the previous case, the existence of such a $\bar s<s_1<s_o$ follows.

\v \n {\bf 6.}   If $\tilde p(s_o)\in D=\{ (p_1,p_2)\,:
~p_1\geq 0\,,~p_2\geq 0\,, ~p_1+p_2\geq 2\kappa_2\}$ for
some $s_o$, then $|\tilde p|\to +\infty$ as $s\to -\infty$.
Indeed, since
$$
{d\over ds} (\tilde p_1+\tilde p_2)= -\tilde p_1^2 -\tilde p_2^2+
\kappa_1 \tilde p_1+\kappa_2 \tilde p_2 \leq
-{1 \over 2} (\tilde p_1+\tilde p_2-2\kappa_2)(\tilde p_1+\tilde p_2)\leq 0
$$
(and the inequality is actually strict when $p_1+p_2= 2\kappa_2$),
we can assume that
there exist $\bar s\leq s_o$ and $\ve>0$ such that
$\tilde p_1(\bar s)+\tilde p_2(\bar s)>2\kappa_2+\ve$. Moreover,
the following holds for any $\sigma<\bar s$:
$$
{d\over ds} (\tilde p_1+\tilde p_2)(\sigma)< -{1 \over 2}
(\tilde p_1 (\bar s)+\tilde p_2(\bar s)-2\kappa_2)
(\tilde p_1(\sigma)+\tilde p_2(\sigma))< -{\ve \over 2}
(\tilde p_1(\sigma)+\tilde p_2(\sigma))\,.
$$
Hence by integrating we find  $(\tilde p_1+\tilde p_2)(s)\geq
\eta e^{-{\ve \over 2}
s}$ for $s<\bar s$ and $\eta>0$. Therefore
$(\tilde p_1+\tilde p_2)\to +\infty$ as $s\to -\infty$.

\v \n {\bf 7.}   If $\tilde p(s_o)\in E=\{ (p_1,p_2)\neq
(0,0)\,: ~p_1\geq 0\,,~p_2\geq 0\,, ~p_1+p_2\leq
2\kappa_1\}$ for some $s_o$, then
from
$$
{d\over ds} (\tilde p_1+\tilde p_2)= -\tilde p_1^2-\tilde p_2^2+
\kappa_1 \tilde p_1+\kappa_2 \tilde p_2 \geq
 -{1 \over 2} (\tilde p_1+\tilde p_2-2\kappa_1)(\tilde p_1+\tilde p_2)\geq 0\,,
$$
it follows, as above, that
either $\tilde p\to 0$ for $s\to -\infty$ or there
exists $\bar s\leq s_o$ such that $\tilde p(\bar s)$ satisfies one of
the previous cases {\bf 3}-{\bf 4}-{\bf 5}.

\v \n {\bf 8.}   If $\tilde p(s_o)\in F=\{ (p_1,p_2)\,:
~p_1\geq 0\,,~p_2\geq 0\,, ~2\kappa_1<p_1+p_2<2\kappa_2\}$
for some $s_o$ and $p\neq\tilde p$, then there exists a small circle $V$
(say with radius smaller than $|\tilde p (s_o)- p(s_o)|$) around the
stable focus
$p\equiv (\kappa_1,\kappa_2)$  such that $\tilde p\notin V$ for $s<s_o$. But
then, looking at the signs of the derivatives of $\tilde p_i$,
as $s\to -\infty$ our solution $\tilde p$ must go away from the whole region
$F$ and there exists $\bar s< s_o$ such that $\tilde p (\bar s)$
is in one of the previous cases.

\midinsert
\vskip 10pt
\centerline{\hbox{\psfig{figure=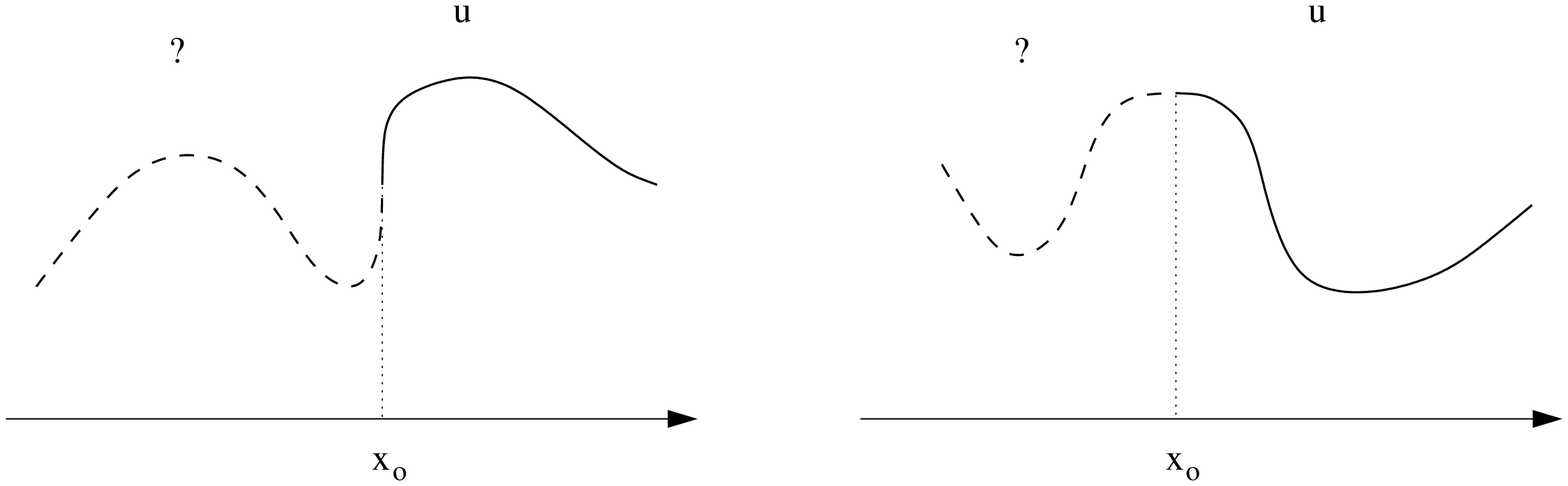 ,width=12cm}}}
\centerline{\hbox{Figure 3a~~~~~~~~~~~~~~~~~~~~~~~~~~~~~~~~~~~~~~~Figure 3b}}
\vskip 10pt
\endinsert

\v \n {\bf 9.}   We can now provide more accurate estimates
on blow-up. Indeed by previous analysis, blow-up of $|\tilde p|$ can
only occur when either $\tilde p_i\to -\infty$ as $s\to +\infty$
 or $\tilde p_i\to +\infty$ as $s\to -\infty$, for some index $i\in\{1,2\}$.
 To fix the ideas, assume $|\tilde p_1|\to\infty$. Then for $s$ sufficiently
large
$$
-{{\tilde p_1^2} \over 2}+(\kappa_1-\kappa_2) \tilde p_1+
\kappa_1 \tilde p_2 <0\,.
$$
Integrating  the inequality ${d\over ds} \tilde p_1< -{{\tilde p_1^2} \over 2} $,
one can conclude that $|\tilde p|\to\infty$ as $s\to s_o$, for some finite
$s_o$ in both cases.
In particular for this $s_o\in\R$ (and $\eta>0$) $\tilde p$ satisfies
$|\tilde p(s)|\geq{\eta \over {|s-s_o|}}$.

In terms of the original variable $x$, one may guess that the corresponding
function $\tilde p=\tilde p(x)$ could be as in Figure~3a and that $u$
may be continued beyond the point where $\tilde p$ blows-up
(say $x_o=x(s_o)$). But this is not the case since such a trajectory
yields a solution defined on the whole real line. Indeed by (4.7)
$${\left|dx\over ds\right|}= \Delta\big(\tilde p(s)\big) \geq
{c_o\over (s_o-s)^2}\,.\eqno(4.17)$$
for some $c_o>0$, and therefore either $x(s)\to +\infty$ as $s\to s_o-$ or
$x(s)\to -\infty$ as $s\to s_o+$.
Therefore, the solution $\tilde u(x)$, corresponding to $\tilde p(x)$,
violates the growth assumptions (2.15) and is not admissible.

\v \n {\bf 10.}   We remark that in case {\bf 7}, the solution
$\tilde p$ can tend to $0$ as $s\to -\infty$. But then for some $c_o>0$
$$
|\tilde p|\leq \tilde p_1+\tilde p_2\leq e^{c_o s}\,.
$$
\n Recalling (4.7) we obtain, in terms of the variable $x$,
$$
\left|{dx\over ds}\right|= \Delta\big(\tilde p(s)\big)=\O(1)\cdot
e^{2 c_o s}\,, \eqno(4.18)$$
$$\lim_{s\to -\infty}x(s)=x_o<\infty\,,$$
for some $x_o\in\R$.
Therefore, to the entire trajectory $s\mapsto \tilde p(s)$,
there corresponds only a portion of the trajectory $x\mapsto \tilde p(x)$,
namely for $x>x_o$.

To prolong the solution $\tilde u$ for $x<x_o$,
we need to construct another trajectory $s\mapsto p(s)$ such that
$\lim_{s\to +\infty} p(s) =0$.
But this trajectory, by previous analysis, will be unbounded,
hence the corresponding $\tilde u (x)$,
will not be admissible.
\v\n {\bf 11.} Next, we consider the case where $\tilde p (s)$
is a discontinuous solution
with admissible jumps. In this case, first of all we can say that $\tilde p$
has no more than 2 jumps. Indeed the set
$\Xi_1=\{ (p_1,p_2)\,: ~p_2< 0\,, ~p_1+p_2\leq 0\}$ is positively invariant
and $\Xi_2=\{ (p_1,p_2)\,: ~p_2< 0\,, ~p_1+p_2> 0\}$ is
negatively invariant. Hence if a jump occurs at $s_o$, either
$\tilde p(s_o+)\in \Xi_1$ or $\tilde p(s_o-)\in \Xi_2$. In the former
case $\tilde p(s)$ has no jumps for $s>s_o$; in the
latter case $\tilde p$ has no jumps for $s<s_o$. This means that there
could be at most two jumps when there exist $s_1<s_2\leq s_3$ such that

\i{$\bullet$} a first jump occurs at $s_1$ and  $\tilde p(s_1-)\in \Xi_2$,
\i{$\bullet$} $\tilde p$ crosses the line $p_1+p_2=0$ at $s_2$,
\i{$\bullet$} a last jump occurs at $s_3$ and $\tilde p(s_3+)\in \Xi_1$.

\n In any case, the corresponding solution $\tilde u$
does not
satisfy (2.15) and is not admissible.  Indeed, we can have only
three situations for a $\tilde p$ with an admissible jump at $s_o$:

\i{(a)} if $\tilde p(s_o-)\in\Xi_2$, then $|\tilde p|\to\infty$ as
$s\to-\infty$;
\i{(b)} if $\tilde p(s_o+)\in\Xi_1$ and $\tilde p_1(s_o+)>0$, then either
$\tilde p(s)$ is continuous for $s<s_o$ (and therefore $|\tilde p|\to\infty$
as $s\to-\infty$) or $\tilde p$ has another jump at $\bar s$ such that
$\tilde p(\bar s-)\in\Xi_2$ (and therefore again $|\tilde p|\to\infty$ as
$s\to-\infty$);
\i{(c)} if $\tilde p(s_o+)\in\Xi_1$ and $\tilde p_1(s_o+)\leq 0$, then
$|\tilde p|\to\infty$ as $s\to +\infty$.

\vs
\n {\it Second Step.} We now extend the proof, in the presence of a sufficiently
small perturbation.
By (4.16), there exist constants $\kappa_1,\kappa_2>0$ such that
$$\big|h_1'(x)-\kappa_1\big|\leq\delta\,,\qquad\qquad
\big|h_2'(x)-\kappa_2\big|\leq\delta\qquad\qquad \hbox{for all }~
x\in\R\,.\eqno(4.19)$$

Let $u(\cdot)$ be the solution constructed in Theorem~2, and let
$\tilde u$ be any other smooth solution of (4.8). Call $p=u'$,
$\tilde p=\tilde u'$ the corresponding gradients, rescaled as before,
and let $V$ be a small open bounded set containing the whole image of
$p$ and the point $(\kappa_1,\kappa_2)$. Of course it is not restrictive
to consider $V$ as circular, say with radius $\rho>0$.

Now we split the proof in three cases.

\v \n CASE 1: $\tilde p(s)\in V$ for every $s$. In this case we
look at the difference $w(s)=\tilde p(s)-p(s)$. We can write
a linear evolution equation for $w$:
$${dw\over ds}=A(s)\,w(s)\,,\eqno(4.20)$$
where the matrix $A$ is the ``average'' matrix

$$
A(s)=\int_0^1 Df(\theta p(s)+(1-\theta) \tilde p(s)) d\theta\,, \eqno(4.21)
$$

\n and $f$ is the vector field at (4.12).

Since $p,\tilde p\in V$, the matrix $A(s)$ is very close to
the Jacobian matrix $Df(\kappa_1,\kappa_2)$, therefore
$${d\over ds}\big|w(s)\big|\leq - K\,\big|w(s)\big|\,,\eqno(4.22)$$
for some constant $K>0$. Indeed $Df(\kappa_1,\kappa_2)$ is negative
definite and, provided $\delta$ (and then $\rho$) is small enough,
$A(s)$ is negative definite too. Hence
$$
2\big|w(s)\big|{d\over ds}\big|w(s)\big|= {d\over
ds}\big|w(s)\big|^2= 2{d\over ds} w(s)\cdot w(s)=2 A(s)w(s)\cdot
w(s)\leq -2 K \big|w(s)\big|^2\,.
$$

Now integrating (4.22), we have for $s<0$,
$$2\rho\geq\big|w(s)\big|\geq e^{-Ks}\big|w(0)\big|\,$$
and, letting $s\to -\infty$, find
$$\big|w(0)\big|\leq \lim_{s\to -\infty}
 \big|w(s)\big| e^{Ks}\leq\lim_{s\to -\infty}2\rho e^{Ks}=0\,.\eqno(4.23)$$
This implies $p(0)=\tilde p(0)$, hence $p=\tilde p$ by the
uniqueness of the Cauchy problem.

\v \n CASE 2: $\tilde p(s_o)\notin V$ for some $s_o$ and, in particular,
$\tilde p(s_o)$ in a small neighbourhood $W$ of the origin. Consider
the linearized system near $(0,0)$

$$\pmatrix{ p_1' \cr p_2'\cr} = H \cdot\pmatrix{ p_1 \cr p_2\cr}
\,\,, \qquad\qquad H=\pmatrix{ {h_1'-h_2'} & h_1' \cr h_2' & {h_2'-h_1'} \cr}\,,$$

\n and notice that the origin is a saddle point for this system.
Indeed $H$ has eigenvalues $\lambda_1,\,\lambda_2$ such that

$$ 0<\sqrt{{3\over {4C^2}}} \leq|\lambda_i|=
\sqrt{(h_1')^2+(h_2')^2-h_1'h_2'}\leq
\sqrt{2C^2-{1\over C^2}}\,, \eqno(4.24)$$

\n where $C$ is the constant in (4.11).
Moreover its eigenvectors $v_1,\,v_2$ form angles
$\alpha_1,\,\alpha_2$ with the positive direction of the $p_1$-axis such that

$$
0<{1 \over C}
\left(\sqrt{d^2+{1 \over {C^2}}}-d\right)\leq
| \tg \alpha_i | = \left|{{\lambda_i+(h_2'-h_1')}\over h_1'}\right| \leq C
\left(\sqrt{d^2+C^2}+d
\right)\,,\eqno(4.25)
$$

\n where $d=(C-{1 \over C})>0$ and $C$ is again from (4.11).

\n Hence, exactly as one can do with saddle points in the autonomous case,
we can prove that there exist four sectors $S_i$, $i=1,\ldots,4$
(see Figure~4), where the following facts hold:

\midinsert \vskip 10pt
\centerline{\hbox{\psfig{figure=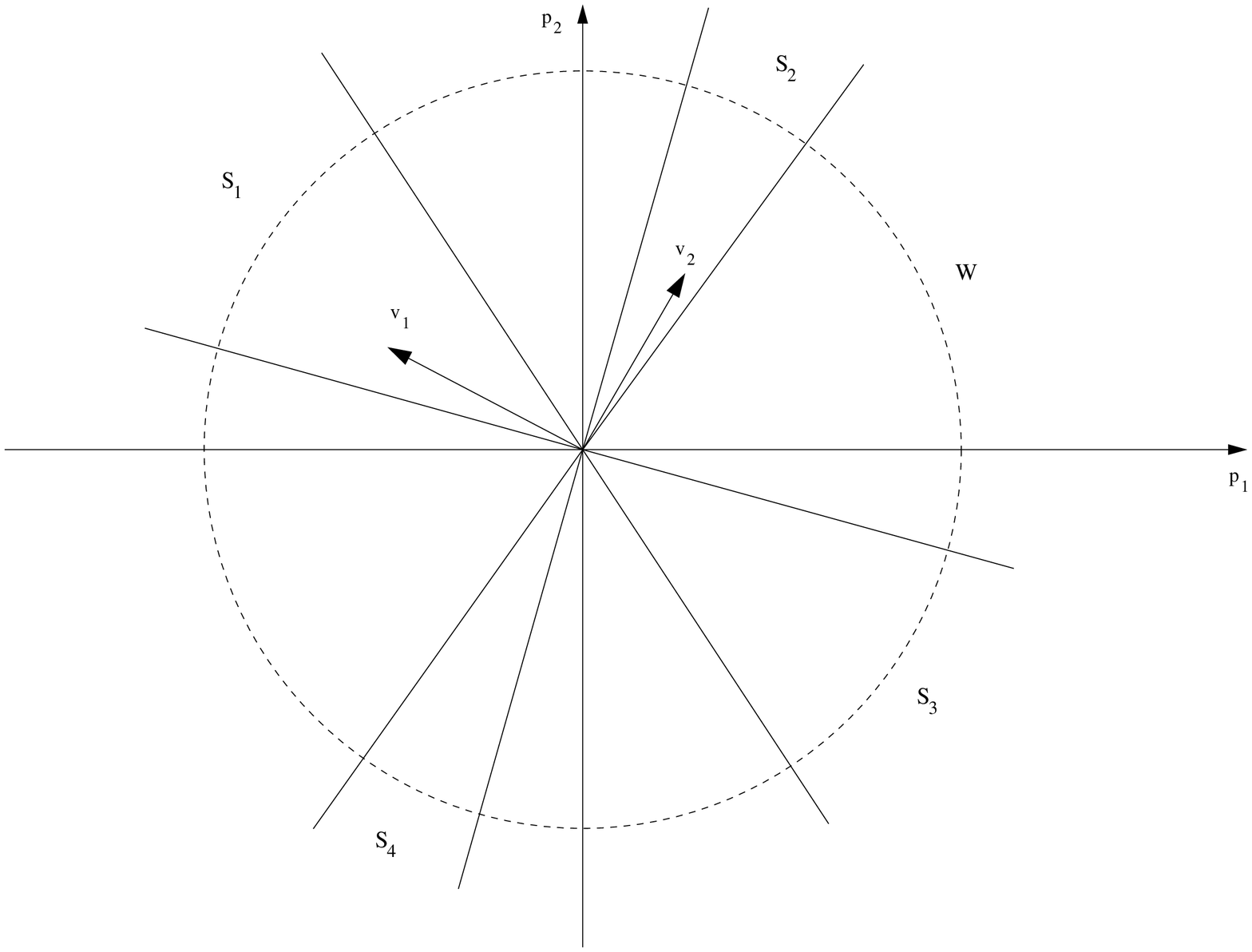,width=10cm}}}
\centerline{\hbox{Figure 4}} \vskip 10pt
\endinsert

\v \n {\bf (a)}   If $\tilde p(s_o)$ is in $S_1$ or $S_3$,
then $|\tilde p(s)|$ grows for $s<s_o$ and the solution moves away
from $W$;
\v \n {\bf (b)}   Both boundaries of $S_2$ and $S_4$ allow
orbits to only exit from those sectors for $s<s_o$;
\v \n {\bf (c)}   If $\tilde p(s_o)\notin S_i$ for all $i=1,\ldots,4$,
then for $s<s_o$ the angle between the vector $(\tilde p_1, \tilde p_2)$
and the $p_1$-axis is strictly monotone, forcing the solution either to
reach $S_1$ or $S_3$, or to move away from $W$;
\v \n {\bf (d)}   Finally, if $\tilde p(s_o)$ is in $S_2$ or $S_4$,
then for $s<s_o$ the solution can tend to the origin. But, since

$$
{d\over ds} (\tilde p_1+\tilde p_2)= -\tilde p_1^2-\tilde p_2^2+
h_1' \tilde p_1+h_2' \tilde p_2 \geq
 -{1 \over 2} (\tilde p_1+\tilde p_2-2C)(\tilde p_1+\tilde p_2)> 0\,,
$$

\n as in the constant case, one obtains an estimate of exponential type of the
decay of $|\tilde p|$.

\v \n CASE 3: $\tilde p(s_o)\notin V$ and $\tilde p(s_o)$ not
in a neighbourhood of the origin. In this case, combining (4.19)
and the continuous dependence of solutions with the estimates of
the constant case (indeed, using (4.11), they remain true),
we can prove that $|\tilde p|\to\infty$ for finite $s$ and
that the rate of blow-up of $|\tilde p|$ can be estimated in the same way
we did in the case of $h_i'\equiv\kappa_i$.

\v \n In any case either $\tilde u\equiv u$ or, in the original coordinates
$x$, $\tilde u$ fails to satisfy (2.15).

\v It remains to prove what happens if $\tilde u$ is an admissible solution
with discontinuous (rescaled) gradient $\tilde p (s)$. Then assume $\tilde p$
has an admissible jump at $s_o$. Using (4.19) it holds, for $\tilde p_1>0$,

$$
\left.{d\over ds} (\tilde p_1+\tilde p_2)\right|_{\tilde p_1+\tilde p_2=0}=
-2\tilde p_1^2+(h'_1-h'_2)\tilde p_1<
-2\tilde p_1^2+(\kappa_1-\kappa_2+2\delta)\tilde p_1
$$

\n and hence, provided $\delta$ small enough, the region $\Xi_1$
(resp. $\Xi_2$) defined in the {\it First Step} is positively (resp.
negatively) invariant also in this setting. Then conclusions made in the
constant case still hold and $\tilde u$ corresponding to $\tilde p$ is not
admissible, since it violates (2.15).
\endproof

\vsk

\n{\medbf 5 - Players with conflicting interests}
\v
We consider here a game for two players, with dynamics (4.1)
and cost functionals as in (4.2).
Contrary to the previous section, we now assume that the
player have conflicting interest. Namely, their running costs
$h_i$ satisfy
$$h'_1(x)\leq 0\leq h_2'(x).\eqno(5.1)$$
We begin with an example showing that in this case
the H-J system can have infinitely many admissible solutions.
Each of these determines a different Nash equilibrium solution
to the differential game.
\v
\n{\bf Example 2.}  Consider the game (4.1)-(4.2),
with
$$h_1(x)= -\kappa x\,,\qquad\qquad h_2(x)=\kappa x\,,\eqno(5.2)$$
for some constant $\kappa>0$ (see Figure~5).

\midinsert \vskip 10pt
\centerline{\hbox{\psfig{figure=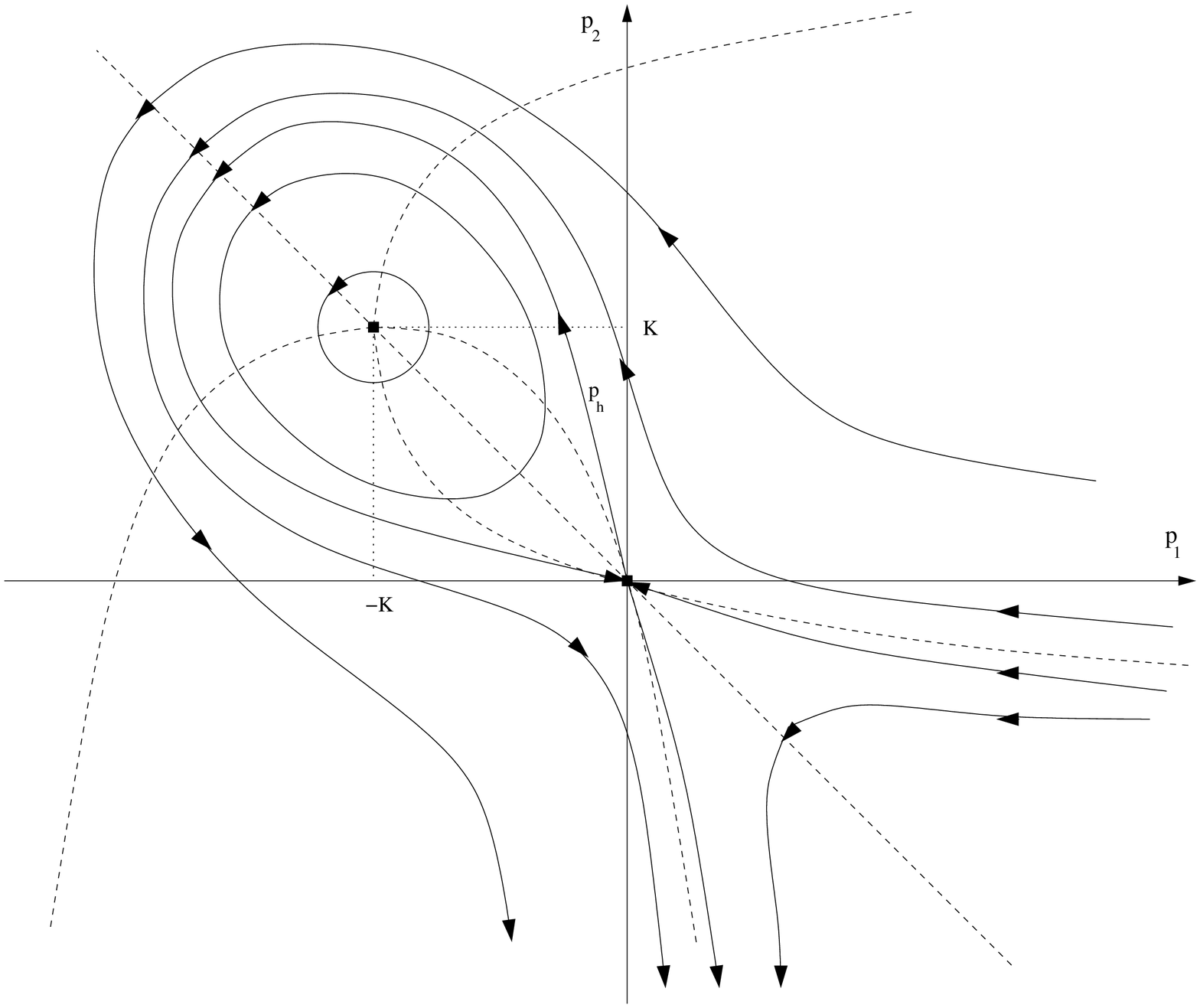,width=12cm}}}
\centerline{\hbox{Figure 5}} \vskip 10pt
\endinsert

In this special case,
the equations (4.8) reduce to
$$\left\{
\eqalign{p_1'
&=-2\kappa p_1-\kappa p_2 -p_1^2
\,,\cr
p_2'
&=~~\kappa p_1+2\kappa p_2 -p_2^2
\,.\cr}\right.\eqno(5.3)$$
The point $\ov P\doteq (-\kappa, \kappa)$ is stationary for the
flow of (5.3). Settinq $q_1\doteq p_1+\kappa$,
$q_2\doteq p_2-\kappa$,
the local behavior of the system near $\ov P$ is described by
$$\left\{
\eqalign{q_1'
&=-\kappa q_2 -q_1^2
\,,\cr
q_2'
&=~\kappa q_1 -q_2^2
\,.\cr}
\right.\eqno(5.4)$$
Notice that
$${dp_2\over dp_1}={dq_2\over dq_1}=0\qquad\hbox{if}\qquad
q_1={q_2^2\over\kappa}\,,$$
$${dp_1\over dp_2}={dq_1\over dq_2}=0\qquad\hbox{if}\qquad
q_2=-{q_1^2\over\kappa}\,,$$
$${dp_1\over dp_2}={dq_1\over dq_2}=1\qquad\hbox{if}\qquad
p_1=-p_2\,.$$
By symmetry across the line $p_1+p_2=0$,
any trajectory passing through a point
$P_\alpha\doteq (-\alpha,\alpha)$ with $0<\alpha<\kappa$
is a closed orbit.
We thus have infinitely many solutions of the H-J equations (5.3),
having bounded, periodic gradients.  Therefore, all of these solutions
are globally Lipschitz continuous and satisfy the growth condition
(2.15).
Notice that the homoclinic orbit $p_h(\cdot)$ starting and ending
at the origin also yields a periodic solution to
the original equation (4.6).
Indeed, to a solution
$p=p(s)$ of (5.3) with
$$\lim_{s\to -\infty}p(s)=\lim_{s\to +\infty}p(s)=0\,,$$
through the reparametrization $x=x(s)$ there corresponds
a solution $p=p(x)$ defined on some bounded interval
$\,]\ell_o, \ell_1[\,$.  This yields a
periodic solution $p=p(x)$ with period $\ell=\ell_1-\ell_o$.

\vs

The main result of this section is concerned with
the existence and uniqueness of admissible solutions.
\v
\n{\bf Theorem 4.} {\it Let any two constants $\kappa_1,\kappa_2$
be given, with
$$\kappa_1<0<\kappa_2\,,\qquad\qquad \kappa_1+\kappa_2\not= 0\,.
\eqno(5.5)$$
Then there exists $\delta>0$ such that the following holds.
If  $\,h_1,h_2$ are smooth functions whose derivatives satisfy
$$\big|h'_1(x)-\kappa_1\big|\leq\delta\,,\qquad\qquad
\big|h'_2(x)-\kappa_2\big|\leq\delta\,,\eqno(5.6)$$ for all $x\in\R$,
then the system of H-J equations (4.3) has a unique admissible
solution.}
\v
\n{\bf Proof.}
We will first consider
the linear case, where $h_i'\equiv \kappa_i$ is constant. Then
we recover the more general case by a perturbation argument.

\v\n   {\bf Existence.}  Assume that
$h_i(x)=\kappa_i\,x$ with
$\kappa_1+\kappa_2>0$, which is not rectrictive. The existence of an
admissible solution for (4.8) is trivial, since we have the constant solution
$p\equiv(\kappa_1,\kappa_2)$, which corresponds to
$$
\big(u_1(x),\, u_2(x)\big)=(\kappa_1 x +\kappa_1\kappa_2 +{\kappa_1^2 \over 2}\,
,~\kappa_2 x +\kappa_1\kappa_2 +{\kappa_2^2 \over 2})\,.\eqno(5.7)
$$

\v Consider now the case of $h'_1,h'_2$ small perturbations of the constants
$\kappa_1,\kappa_2$. Notice that, in the previous case,
every ball $B(\kappa, R)$ around $\kappa=(\kappa_1,\kappa_2)$ with radius
$R< {\sqrt{2} \over 2}(\kappa_1+\kappa_2)$ was positively invariant for
the flow of (4.8).

\n Indeed, setting $q_i= p_i-\kappa_i$, the system becomes
$$
\left\{\eqalign{
q_1'&=-(\kappa_1+\kappa_2) q_1+\kappa_1 q_2 -q_1^2\,,\cr
q_2'&=\kappa_2 q_1-(\kappa_1+\kappa_2) q_2 -q_2^2\,.\cr
}\right.\eqno(5.8)
$$

\n and it holds

$$
{d\over ds} {|q|^2 \over 2} =
- q_1^3 - q_2^3 - (\kappa_1+\kappa_2) (q_1^2+q_2^2-q_1q_2)=
- (q_1^2+q_2^2-q_1q_2)(\kappa_1+\kappa_2+q_1+q_2)\,.
$$

\n Now, since $|q|\leq R< {\sqrt{2} \over 2}(\kappa_1+\kappa_2)$ ensures
$\kappa_1+\kappa_2+q_1+q_2>0$, one can conclude that

$$
{d\over ds} {|q|^2 \over 2} =
- (q_1^2+q_2^2-q_1q_2)(\kappa_1+\kappa_2+q_1+q_2)
\leq - {|q|^2 \over 2}(\kappa_1+\kappa_2+q_1+q_2) <0\,, \eqno(5.9)
$$

\n and this prove the positively invariance of such a ball $B$.

Then, provided $\delta$ is small enough, we can choose one of these
balls as
a neighborhood $U$ of $(\kappa_1,\kappa_2)$ positively invariant also for
the perturbed system (i.e. $h'_i \not\equiv \kappa_i$). Once we found such a
compact, positively invariant set $\overline U$, we can repeat the existence proof of
Theorem~2:

\i{a.} Consider $p^{(\nu)}\colon \,[-\nu,\infty[\,\to\R^2 $
solution of the Cauchy problem with initial datum
$p^{(\nu)}(-\nu)=(\kappa_1,\kappa_2)$;

\i{b.} By positive invariance, $p^{(\nu)}(x)\in U$ for $x>-\nu$.
We then extend the function $p^{(\nu)}$ to
the whole real line by setting $p^{(\nu)}(x)\equiv (\kappa_1,\kappa_2)$
for $x<-\nu$;

\i{c.} By uniform boundedness and equicontinuity, the sequence
$p^{(\nu)}$ admits a subsequence converging to a uniformly continuous
function $p:\R\mapsto U$.  Clearly, this limit function
$p(\cdot)$ provides a global, bounded solution to the system (4.8). In turn, this
yields an admissible solution $u(\cdot)$ to (4.6).

\v\n   {\bf Uniqueness.} {\it First Step.} Let
$h'_i\equiv \kappa_i$ and $\kappa_1+\kappa_2>0$. In order to prove that
the previously found solution is the only one that satisfies (A1)-(A3),
we assume that $\tilde u$ is another solution of the system (4.3),
whose gradient will be denoted by $\tilde p$. Figure~6 depicts
possible trajectories $s\mapsto \tilde p(s)$ of the planar system (4.8).
We remark that:

\midinsert \vskip 10pt
\centerline{\hbox{\psfig{figure=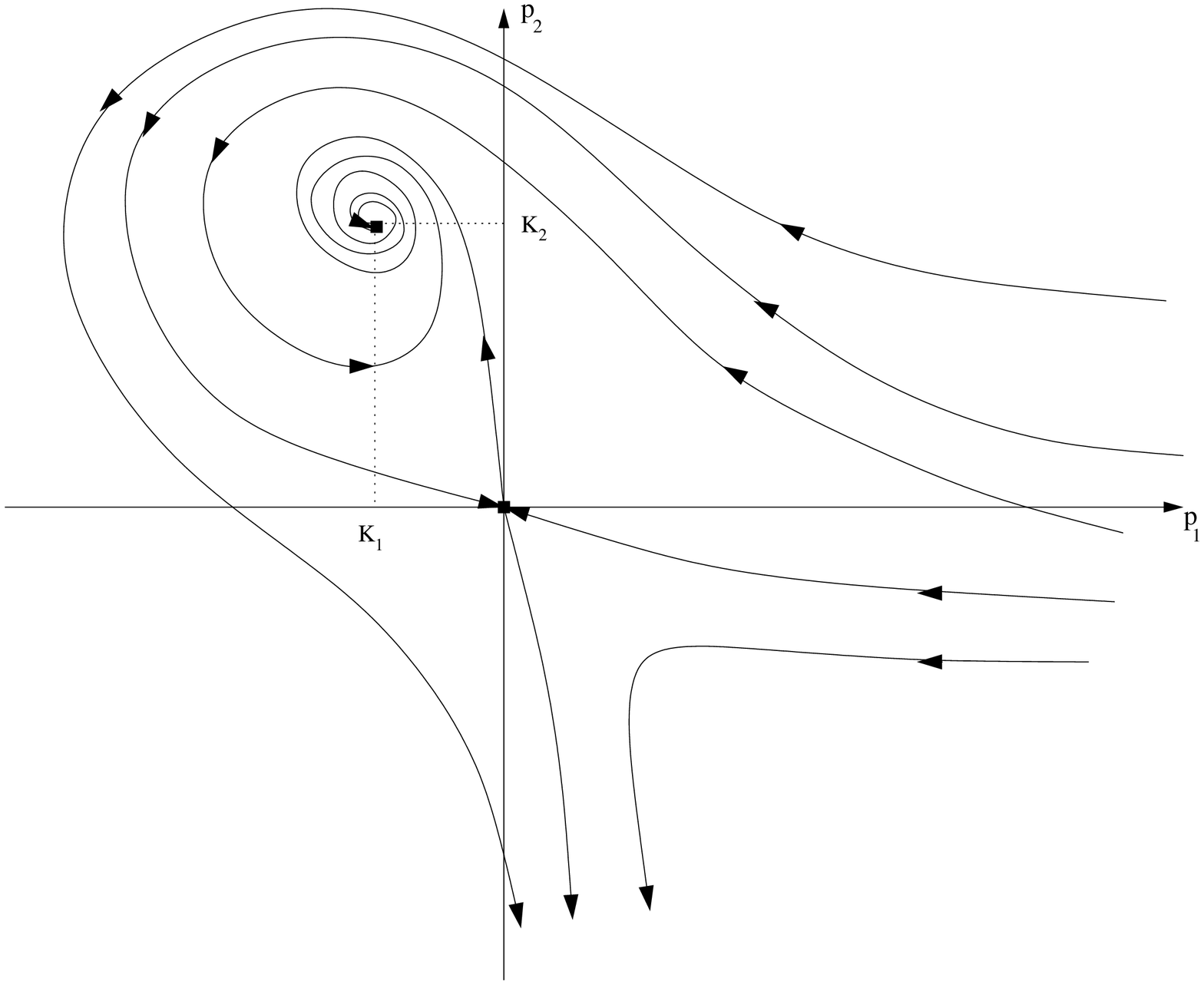,width=12cm}}}
\centerline{\hbox{Figure 6}} \vskip 10pt
\endinsert

\v \n {\bf 1.}   The regions $\{ (p_1,p_2)\neq
(0,0)\,: ~p_1\geq 0\,,~p_2\leq 0\,, ~p_1+p_2\leq 0\}$ and
$\{ (p_1,p_2)\neq (0,0)\,: ~p_1\geq 0\,,~p_2\leq 0\,,
~p_1+p_2\geq 0\}$ are positively and negatively
invariant for the flow of (4.8), respectively.

\v \n {\bf 2.}   If $\tilde p(s_o)\in\{ (p_1,p_2)\neq
(0,0)\,: ~p_1\geq 0\,,~p_2\leq 0\,, ~p_1+p_2\leq 0\}$ for
some $s_o$, then $|\tilde p|\to +\infty$ as $s\to +\infty$.
Indeed, since
$$
{d\over ds} (\tilde p_1+\tilde p_2)= -\tilde p_1^2-
\tilde p_2^2+\kappa_1 \tilde p_1+\kappa_2 \tilde p_2 <
-(\tilde p_1+\tilde p_2)^2 < 0\,,
$$
we can assume there exist $\bar s\geq s_o$ and $\ve>0$ such that
$\tilde p_1(\bar s)+\tilde p_2(\bar s)<-\ve$. Moreover,
for any $\sigma>\bar s$ we have
$$
{d\over ds} (\tilde p_1+\tilde p_2) (\sigma)\leq -(\tilde p_1
(\sigma)+\tilde p_2 (\sigma))^2 < -(\tilde p_1 (\bar s)+ \tilde p_2(\bar
s))(\tilde p_1 (\sigma)+\tilde p_2 (\sigma))< \ve (\tilde p_1
(\sigma)+\tilde p_2 (\sigma))\,.
$$
After an integration, we find $(\tilde p_1+\tilde p_2)(s)\leq -\eta
e^{\ve s}$ for $s>\bar s$ (and $\eta>0$) and hence $(\tilde p_1+\tilde p_2)\to
-\infty$ as $s\to +\infty$.

\v \n {\bf 3.}   If $\tilde p(s_o)\in\{ (p_1,p_2)\neq
(0,0)\,: ~p_1\geq 0\,,~p_2\leq 0\,, ~p_1+p_2\geq 0\}$ for
some $s_o$, then $|\tilde p|\to +\infty$ as $s\to -\infty$.
Indeed, reasoning as above, we can assume there exist $\bar s\leq s_o$
and $\ve>0$ such that $\tilde p_1(\bar s)+\tilde p_2(\bar
s)>\ve$ and the following holds for any $\sigma<\bar s$:
$$
{d\over ds} (\tilde p_1+\tilde p_2) (\sigma)\leq -(\tilde p_1
(\sigma)+\tilde p_2 (\sigma))^2 < -(\tilde p_1 (\bar s)+\tilde p_2 (\bar
s))(\tilde p_1 (\sigma)+\tilde p_2 (\sigma))< -\ve (\tilde p_1
(\sigma)+\tilde p_2(\sigma))\,.
$$
This implies $(\tilde p_1+\tilde p_2)(s)\geq \eta
e^{-\ve s}$ for $s<\bar s$ (and $\eta>0$), hence  $(\tilde p_1+\tilde p_2)\to
+\infty$ as $s\to -\infty$.

\v \n {\bf 4.}   If $\tilde p(s_o)\in\{ (p_1,p_2)\,:
~p_1> 0\,,~p_2> 0\}$ for some $s_o$, then $|\tilde p|\to
+\infty$ as $s\to -\infty$. Indeed, let $\ve>0$ such that
$\tilde p_1(s_o)> \ve$. Since
$$
{d\over ds} \tilde p_1 = -\tilde p_1^2+
(\kappa_1-\kappa_2)\tilde p_1+ \kappa_1 \tilde p_2 <
-\tilde p_1^2< 0\,,
$$
it is sufficient to observe that for $\sigma< s_o$
$$
{d\over ds} \tilde p_1 (\sigma)\leq
-\tilde p_1(s_o)\tilde p_1(\sigma) < -\ve \tilde p_1(\sigma)\,.
$$
Hence, integrating, $\tilde p_1(s)\geq \eta e^{-\ve s}$ for
$s<s_o$ (and $\eta>0$) and either $\tilde p_1\to +\infty$ as $s\to -\infty$ or
there exists $\bar s<s_o$ such that $\tilde p$ is in the previous
case.

\v \n {\bf 5.}   If $\tilde p(s_o)\in\{ (p_1,p_2)\,:
~p_1< 0\,,~p_2< 0\}$ for some $s_o$, then $|\tilde p|\to
+\infty$ as $s\to +\infty$. Here we can repeat the
argument of {\bf 4.} with $\tilde p_2$ in place of $\tilde p_1$.

\v \n {\bf 6.}   Let $\tilde p(s_o)\in\{ (p_1,p_2)\neq
(0,0)\,: ~p_1\leq 0\,,~p_2\geq 0\}$ for some $s_o$ and
set $\hat p$ as the unique solution in this region that tends to
the origin as $s\to +\infty$. Notice that, as $s\to -\infty$, either
$\hat p(s)$ crosses the $p_2$-axis or $\hat p_2\to+\infty$. Then:

\i{$\bullet$} if $\tilde p=\hat p$, then as stated above either
there exists $\bar s<s_o$ such that $\tilde p(\bar s)$ is in the case {\bf 4},
or $\tilde p_2\to\infty$ as $s\to - \infty$. In both cases
$|\tilde p|\to \infty$ as $s\to -\infty$.

\i{$\bullet$} if $\tilde p(s_o)$ belongs to the region between $\hat p$ and the
$p_2$-axis, then there could be only three possibilities: either $\tilde p$
is the unique solution that tends to the origin as $s\to -\infty$, or
$\tilde p_2\to\infty$ as $s\to -\infty$ without $\tilde p$ crosses $p_2$-axis
(and, of course, this can only happen if $\hat p$ does not cross it too),
or there exists $\bar s<s_o$ such that $\tilde p(\bar s)$ is in the case
{\bf 4.} above.  In the former case we will estimate the decay of $|\tilde p|$
in {\bf 8}; in the latter ones $|\tilde p|\to \infty$ as $s\to -\infty$.

\i{$\bullet$} if $\tilde p(s_o)$ doesn't belong to the region between
$\hat p$ and the $p_2$-axis, then either $\tilde p_2\to\infty$ as
$s\to -\infty$ or there exists $\bar s<s_o$ such that
$\tilde p(\bar s)$ is in case {\bf 5} above
(and this is possible only if also $\hat p(s)$ crosses the $p_2$-axis).
In both situations, again, $|\tilde p|\to \infty$ as $s\to -\infty$.

\v \n {\bf 7.}   We can now provide more accurate estimates
on the blow-up rate. Indeed by previous analysis, as in Theorem~3, a blow-up of
$|\tilde p|$ can only occur when either $\tilde p_i\to -\infty$ as
$s\to +\infty$ or $\tilde p_i\to +\infty$ as $s\to -\infty$
(for some index $i$). Hence, exactly as before, we can prove that
there exists $s_o\in\R$ (and $\eta>0$) such that
$|\tilde p(s)|\geq{\eta \over {|s-s_o|}}$.
In terms of the original variable $x$, such a trajectory yields a solution
defined on the whole real line, because by (4.7)
$${\left|dx\over ds\right|}= \Delta\big(\tilde p(s)\big) \geq
{c_o\over (s_o-s)^2}\,.\eqno(5.10)$$
for some $c_o>0$ and therefore either $x(s)\to+\infty$ as $s\to s_o-$ or
$x(s)\to-\infty$ as $s\to s_o+$.
In conclusion, the solution $\tilde u(x)$ which corresponds to  $\tilde p(x)$
violates the growth condition (2.15), and hence it is not admissible.

\v \n {\bf 8.}   Notice that only in case {\bf 6}-(ii), where
$\tilde p$ is the unique solution that tends to $0$ as $s\to -\infty$,
we have a solution that could remain bounded in the whole $\R$. But in this
case, we shall have as $s\to-\infty$
$$
|\tilde p|(s)\leq (\tilde p_2-\tilde p_1) (s)\leq \gamma e^{c_o s}\,,\eqno(5.11)
$$
for some $\gamma, c_o>0$.
Indeed studying the linearized system near the origin we
see that $\tilde p$ tends to $(0,0)$ along the direction
$(1, {{\kappa_2-\kappa_1+\sqrt{\kappa_1^2+\kappa_2^2-\kappa_1\kappa_2}}
\over \kappa_1})$. Then there exists $\bar s$ such that for $s<\bar s$
the following holds:
$$
\tilde p_2(s)> (1+{\sqrt{2}\over 2})
{{\kappa_2-\kappa_1}\over \kappa_1}\,\tilde p_1(s)=
\beta\, {{\kappa_2-\kappa_1+\sqrt{\kappa_1^2+\kappa_2^2-\kappa_1\kappa_2}}
\over \kappa_1}\,\tilde p_1(s)\,,\eqno(5.12)
$$
where
$$
\beta= {{(1+{\sqrt{2}\over 2})
(\kappa_2-\kappa_1)}\over {\kappa_2-\kappa_1+
\sqrt{\kappa_1^2+\kappa_2^2-\kappa_1\kappa_2}}}\,,
\qquad \beta\in (0,1)\,.\eqno(5.13)
$$
Notice that, setting
$$
\eqalign{\alpha &={{\kappa_1-2\kappa_2-(\kappa_2-2\kappa_1)\,\beta\,
{{\kappa_2-\kappa_1+\sqrt{\kappa_1^2+\kappa_2^2-\kappa_1\kappa_2}}
\over\kappa_1}}
\over
{1-\beta\, {{\kappa_2-\kappa_1+\sqrt{\kappa_1^2+\kappa_2^2-\kappa_1\kappa_2}}
\over \kappa_1}}}=\cr
&={{\kappa_1-2\kappa_2-(\kappa_2-2\kappa_1)(1+{\sqrt{2}\over 2}){{
\kappa_2-\kappa_1}\over\kappa_1}}
\over
{1-(1+{\sqrt{2}\over 2}){{\kappa_2-\kappa_1}\over \kappa_1}}} >0\,, \cr }
\eqno(5.14)
$$
we obtain exactly
$$
\beta\,{{\kappa_2-\kappa_1+\sqrt{\kappa_1^2+\kappa_2^2-\kappa_1\kappa_2}}
\over \kappa_1}={{\kappa_1-2\kappa_2-\alpha}\over{\kappa_2-2\kappa_1-\alpha}}\,.\eqno(5.15)
$$
Hence, for $s<\bar s$,
$$\tilde p_2 (s)>{{\kappa_1-2\kappa_2-\alpha}
\over{\kappa_2-2\kappa_1-\alpha}}\, \tilde p_1 (s)\,,\eqno(5.16)$$
i.e. $(\kappa_2-2\kappa_1)\tilde p_2 - (\kappa_1-2\kappa_2) \tilde p_1 >
\alpha (\tilde p_2 - \tilde p_1)$.
Recalling that $|\tilde p|\to 0$ as $s\to -\infty$, which implies
the existence of $c_o>0$ and $\hat s$ such that
$\alpha-\tilde p_1(s)-\tilde p_2(s)>c_o$ for any $s<\hat s$, we find
$$
{d\over ds} (\tilde p_2-\tilde p_1)= \tilde p_1^2-
\tilde p_2^2+(\kappa_2-2\kappa_1)\tilde p_2 - (\kappa_1-2\kappa_2) \tilde p_1 >
(\alpha-\tilde p_1-\tilde p_2)(\tilde p_2 - \tilde p_1)>
c_o(\tilde p_2 - \tilde p_1)\,,
$$
for $s$ small enough (namely $s<\min\{\bar s,\hat s\}$). Integrating we
find  $(\tilde p_2-\tilde p_1) (s)\leq \gamma e^{c_o s}$
($\gamma>0$) and hence (5.11) is proved.
Next, recalling (4.7), in terms of the variable $x$
we obtain
$$
\left|{dx\over ds}\right|= \Delta\big(\tilde p(s)\big)=\O(1)\cdot
e^{2 c_o s}\,, \eqno(5.17)$$
$$\lim_{s\to -\infty}x(s)=x_o<\infty\,,$$
for some $x_o\in\R$.
Therefore, to the entire trajectory $s\mapsto \tilde p(s)$,
there corresponds only a portion of the trajectory $x\mapsto \tilde p(x)$,
namely the values for $x>x_o$.
To extend this trajectory also on the half line $\,]-\infty, x_o]$, we need
to construct
another trajectory $s\mapsto p(s)$ with $\lim_{s\to +\infty} p(s) =0$.
But any such trajectory, by previous analysis, will yield a solution
$\tilde u (x)$,
which violates the sublinear growth condition (2.15) as $x\to -\infty$
and is not admissible.

\vs

\n {\it Second Step.} Next, we prove uniqueness of the admissible
solution the case where $h'_i$ is not constant.  Let
$u(\cdot)$ be the solution constructed before, with $p=u'$
remaining in a small disc $V$, centered at
$(\kappa_1,\kappa_2)$ with radius $\rho>0$, positively invariant
for the flow of (4.8).
Moreover, let
$\tilde u$ be any other smooth solution of (4.3).
We split the proof in three cases.

\v \n CASE 1: $\tilde p(s)\in V$ for every $s$. In this case,
as in Theorem~3, we look at the difference $w(s)=\tilde p(s)-p(s)$ and
at the linear evolution equation for $w$:
$${dw\over ds}=A(s)\,w(s)\,,\eqno(5.18)$$
where $A$ is the averaged matrix
$$
A(s)=\int_0^1 Df(\theta p(s)+(1-\theta) \tilde p(s)) d\theta\,,\eqno(5.19)
$$
and $f$ is the vector field describing our system, as in (4.12).
Since $p,\tilde p\in V$, the matrix $A(s)$ is very close to
the Jacobian matrix $Df(\kappa_1,\kappa_2)$, therefore
$${d\over ds}\big|w(s)\big|\leq - K\,\big|w(s)\big|\,,\eqno(5.20)$$
for some constant $K>0$. Indeed,
$$
Df(\kappa_1,\kappa_2) x\cdot x < -{{\kappa_1+\kappa_2}\over 2}|x|^2
.\eqno(5.21)
$$
Provided that $\delta,\rho>0$ are small enough,
there will exist $K>0$ such that $A(s) x\cdot x < - K|x|^2$.
But then, exactly as in Theorem~3, (5.20) implies $p(0)=\tilde p(0)$ and
hence $p=\tilde p$ by the uniqueness of the Cauchy problem.

\v \n CASE 2: $\tilde p(s_o)\notin V$ for some $s_o$ and, in particular,
$\tilde p(s_o)$ in a small neighborhood $W$ of the origin. Consider
the linearized system
$$\pmatrix{ p_1' \cr p_2'\cr} = H \cdot\pmatrix{ p_1 \cr p_2\cr}
\,\,, \qquad\qquad H=\pmatrix{ {h_1'-h_2'} & h_1' \cr h_2' & {h_2'-h_1'} \cr}\,,$$
and notice that the origin is again a saddle point for this system.
Indeed $H$ has eigenvalues $\lambda_1,\,\lambda_2$ such that, recalling (5.6)
and provided $\delta< {1\over 2} \min \{-\kappa_1,\kappa_2\}$,
$$ 0<{\sqrt{2} \over 2} (\kappa_2-\kappa_1-2\delta)
\leq|\lambda_i|= \sqrt{(h_1')^2+(h_2')^2-h_1'h_2'}\leq
\kappa_2-\kappa_1+2\delta
\,. \eqno(5.22) $$
Moreover its eigenvectors $v_1,\,v_2$ form angles
$\alpha_1,\,\alpha_2$ with the positive direction of the $p_1$-axis such that
$$
0\geq \tg\, \alpha_1 = {{\lambda_1+(h_2'-h_1')}\over h_1'} >
 \tg \,\alpha_2 = {{\lambda_2+(h_2'-h_1')}\over h_1'} \,.
$$
More precisely, for $\delta$ small enough, we have
$$
0\geq \tg \,\alpha_1 > (1-{\sqrt{2}\over 2})\,
{{\kappa_2-\kappa_1+2\delta}\over {\kappa_1+\delta}} >
(1+{\sqrt{2}\over 2})\,
{{\kappa_2-\kappa_1-2\delta}\over {\kappa_1-\delta}}
> \tg\, \alpha_2\,. \eqno(5.23)
$$
Hence, as in the previous proof of
Theorem~3, we show the existence of four sectors
$S_i$, $i=1,\ldots,4$ (see Figure~7), where the following holds:

\midinsert \vskip 10pt
\centerline{\hbox{\psfig{figure=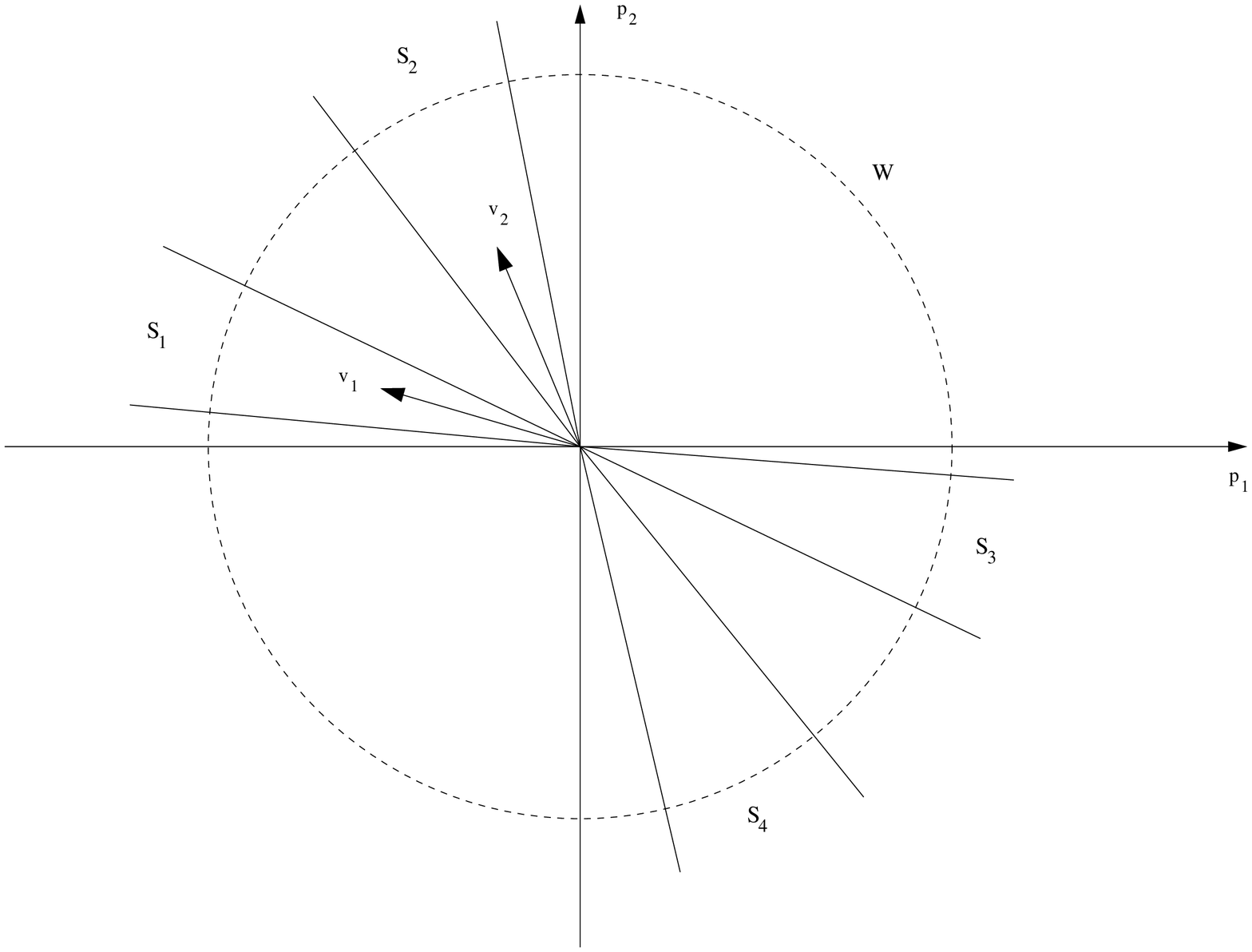,width=10cm}}}
\centerline{\hbox{Figure 7}} \vskip 10pt
\endinsert

\v \n {\bf (a)}   If $\tilde p(s_o)$ is in $S_1$ or $S_3$,
then $|\tilde p(s)|$ grows for $s<s_o$ and the solution moves away
from $W$;
\v \n {\bf (b)}   Both boundaries of $S_2$ and $S_4$ allow
orbits to only exit from those sectors for $s<s_o$;
\v \n {\bf (c)}   If $\tilde p(s_o)\notin S_i$ for all $i=1,\ldots,4$,
then for $s<s_o$ the angle between the vector $(\tilde p_1, \tilde p_2)$
and the $p_1$-axis is strictly monotone, forcing the solution either to
reach $S_1$ or $S_3$, or to move away from $W$;

\v \n {\bf (d)}   Finally, if $\tilde p(s_o)$ is in $S_2$ or $S_4$,
then for $s<s_o$ the solution can tend to the origin. But

$$
\eqalign{{d\over ds} (\tilde p_2-\tilde p_1)&= \tilde p_1^2-\tilde p_2^2+
(h_2'-2h_1') \tilde p_2-(h_1'-2h_2') \tilde p_1 >\cr
&> \tilde p_1^2-\tilde p_2^2+
(\kappa_2-2\kappa_1-3\delta) \tilde p_2-(\kappa_1-2\kappa_2-3\delta)
\tilde p_1\,,\cr}
$$

\n and, provided $\delta$ is small enough, we can use (5.23) to find
$\alpha>0$ such that

$$
{d\over ds} (\tilde p_2-\tilde p_1)> \tilde p_1^2-\tilde p_2^2+
(\kappa_2-2\kappa_1-3\delta) \tilde p_2-(\kappa_1-2\kappa_2-3\delta)
\tilde p_1 > (\alpha -\tilde p_1-\tilde p_2) (\tilde p_2 -\tilde p_1)\,.\eqno(5.24)
$$

\n Hence an estimate of exponential type of the
decay of $|\tilde p|$ follows as in (5.11).

\v \n CASE 3: $\tilde p(s_o)\notin V$ and $\tilde p(s_o)$ not
in a neighbourhood of the origin. In this case, combining (5.6)
and the continuous dependence of solutions with the estimates of
the constant case (indeed they are true also in this more general setting),
we can prove that $|\tilde p|\to\infty$ for finite $s$ and
that the rate of blow-up of $|\tilde p|$ can be estimated in the same way
we did in the {\it First Step}.

\v \n In any case either $\tilde u\equiv u$ or, in the original coordinates
$x$, $\tilde u$ fails to satisfy (2.15).

\v Finally, we rule out the possibility that the gradient
$\tilde p=\tilde u'$ has jumps.
Looking at the phase portrait in Figure 6, we see that
after a one or at most two
admissible jump, the values of $\tilde p$ must fall within
the positively invariant region
$\{ (p_1,p_2)\,: ~p_2< 0\,, ~p_1+p_2< 0\}$. It follows that $\tilde p$ cannot
have any more jumps, and the estimates in {\bf 2},
{\bf 5} and {\bf 7} (together with their analogs in the non-constant case)
imply that $\tilde u (x)$ violates (2.15) as $x\to +\infty$.
Therefore, $\tilde u$ is not an admissible solution.
\endproof

\vsk

\c{\medbf References}
\v
\item{[AFJ]}
H.~Abou-Kandil, G.~Freiling and G.~Jank,
Solution and asymptotic behavior of coupled Riccati equations
in jump linear systems,
{\it IEEE Trans. Automat. Control} {\bf 39}  (1994), 1631--1636.
\v
\item{[BC]}
{M.~Bardi and I.~Capuzzo-Dolcetta},
{\it Optimal control and viscosity solutions of Hamilton-Jacobi-Bellman
equations},
Birkh\"auser, Boston, 1997.
\v
\i{[BS1]} A.~Bressan and W.~Shen,
Small BV solutions of hyperbolic non-cooperative differential
games, {\it SIAM J. Control Optim.} {\bf 43} (2004), 104-215
\v
\i{[BS2]} A.~Bressan and W.~Shen,
 Semi-cooperative strategies for differential games,
{\it Intern. J. Game Theory} {\bf 32} (2004), 561-593.
\v
\item{[EW]}
{J.~C.~Engwerda and A.~J.~T.~M.~Weeren},
{The open-loop Nash equilibrium in LQ-games revisited},
{\it CentER Discussion Paper} {\bf 9551}, Tilburg University, 1995.
\v
\item{[F]}
A.~Friedman,
{\it Differential Games},
Wiley-Interscience, 1971.
\v
\item{[PMC]}
{G.~P.~Papavassilopoulos, J.~V.~Medani\'c and J.~B.~Cruz,~Jr.},
{On the existence of Nash strategies and solutions to coupled
Riccati equations in linear-quadratic games},
{\it  J. Optim. Theory Appl.} {\bf 28}  (1979), 49--76.
\v
\item{[WSE]}
{A.~J.~T.~M.~Weeren, J.~M.~Schumacher and J.~C.~Engwerda},
{Asymptotic analysis of linear feedback Nash equilibria in
nonzero-sum linear-quadratic differential games},
{\it  J. Optim. Theory Appl.} {\bf 101}  (1999), 693--722.
\bye